\newif\ifpdf
\ifx\pdfoutput\undefined
\pdffalse % we are not running PDFLaTeX
\else
\pdfoutput=1 % we are running PDFLaTeX
\pdftrue \fi

\newif\iffinal
\finaltrue % Final version

\documentclass[reqno,twoside,11pt]{amsart}
\iffinal\else\usepackage[notref,notcite]{showkeys}\fi
\usepackage{amsmath}
\usepackage{amsfonts}
\usepackage{amssymb}
\usepackage{verbatim}
\usepackage{epsfig}
\usepackage{setspace}

\IfFileExists{epsf.def}{\input epsf.def}{\usepackage{epsf}}

\IfFileExists{myowntimes.sty}{\usepackage{myowntimes}\usepackage{mathrsfs}}
	{\usepackage{mathpazo}\usepackage{mathrsfs}}

\DeclareFontFamily{OT1}{eusb}{} \DeclareFontShape{OT1}{eusb}{m}{n} {<5> <6> <7> <8> <9> <10> <11> <12> <14.4> eusb10}{}
\DeclareMathAlphabet{\eusb}{OT1}{eusb}{m}{n}

\DeclareFontFamily{OT1}{eusm}{} \DeclareFontShape{OT1}{eusm}{m}{n} {<5> <6> <7> <8> <9> <10> <11> <12> <14.4> eusm10}{}
\DeclareMathAlphabet{\eusm}{OT1}{eusm}{m}{n}

\DeclareFontFamily{OT1}{eufm}{} \DeclareFontShape{OT1}{eufm}{m}{n} {<5> <6> <7> <8> <9> <10> <11> <12> <14.4> eufm10}{}
\DeclareMathAlphabet{\mathfrak}{OT1}{eufm}{m}{n}

\DeclareFontFamily{OT1}{fraktura}{}
\DeclareFontShape{OT1}{fraktura}{m}{n} {<5> <6> <7> <8> <9> <10> <11> <12> <13> <14.4> [1.1] eufm10}{}
\DeclareMathAlphabet{\fraktura}{OT1}{fraktura}{m}{n}

\DeclareFontFamily{OT1}{cmfi}{} \DeclareFontShape{OT1}{cmfi}{m}{n} {<5> <6> <7> <8> <9> <10> <11> <12> <13> <14.4> [0.9] cmfi10}{}
\DeclareMathAlphabet{\cmfi}{OT1}{cmfi}{b}{n}

\DeclareFontFamily{OT1}{cmss}{} \DeclareFontShape{OT1}{cmss}{m}{n} {<5> <6> <7> <8> <9> <10> <11> <12> <13> <14.4> cmss10}{}
\DeclareMathAlphabet{\cmss}{OT1}{cmss}{m}{n}
\newtheoremstyle{thm}{1.5ex}{1.5ex}{\itshape\rmfamily}{} {\bfseries\rmfamily}{}{2ex}{}

\newtheoremstyle{def}{1.5ex}{1.5ex}{\slshape\rmfamily}{} {\bfseries\rmfamily}{}{2ex}{}

\newtheoremstyle{rem}{1.3ex}{1.3ex}{\rmfamily}{} {\itshape}
{} {1.5ex}{}

\theoremstyle{thm}
\newtheorem{theorem}{Theorem}[section]
\newtheorem{lemma}[theorem]{Lemma}
\newtheorem{proposition}[theorem]{Proposition}

\newtheorem*{Main Theorem}{Main Theorem.}
\newtheorem{corollary}[theorem]{Corollary}
\newtheorem*{special theorem}{Lindeberg-Feller Theorem for Martingales}

\theoremstyle{def}

\theoremstyle{rem}
\newtheorem{remark}{{\itshape Remark}}[]

\numberwithin{equation}{section}

\renewcommand{\section}{\secdef\sct\sect}
\newcommand{\sct}[2][default]{%
\refstepcounter{section}
\addcontentsline{toc}{section}{{\tocsection {}{\thesection}{\!\!\!\!#1\dotfill}}{}}
\vspace{0.7cm}
\centerline{\scshape\thesection.\ #1} \nopagebreak \vspace{0.2cm}}
\newcommand{\sect}[1]{%
\vspace{0.4cm} \centerline{\large\scshape\rmfamily #1}
\vspace{0.2cm}}

\renewcommand{\subsection}{\secdef\subsct\sbsect}
\newcommand{\subsct}[2][default]{\refstepcounter{subsection}
\addcontentsline{toc}{subsection}
{{\tocsection{\!\!}{\hspace{1.2em}\thesubsection}{\!\!\!\!#1\dotfill}}{}}
\nopagebreak\vspace{0.45\baselineskip} {\flushleft\bf
\thesubsection~\bf #1.~}
\\*[3mm]\noindent
\nopagebreak}
\newcommand{\sbsect}[1]{\vspace{0.1cm}\noindent
\textbf{#1.~}\vspace{0.1cm}}

\renewcommand{\subsubsection}{%
\secdef \subsubsect\sbsbsect}
\newcommand{\subsubsect}[2][default]{%
\refstepcounter{subsubsection} 
\addcontentsline{toc}{subsubsection}{{\tocsection{\!\!}
{\hspace{3.05em}\thesubsubsection}{\!\!\!\!#1\dotfill}}{}}
\nopagebreak
\vspace{0.15\baselineskip} \nopagebreak {\flushleft\rmfamily
\itshape\thesubsubsection
\ \rmfamily #1\/.}\ }
\newcommand{\sbsbsect}[1]{\vspace{0.1cm}\noindent
\rmfamily \itshape
\arabic{section}.\arabic{subsection}.\arabic{subsubsection} \
\sffamily #1\/.\ }

\iffinal

\else

\fi

\renewcommand{\caption}[1]{%
\vglue0.5cm
\refstepcounter{figure}
\begin{minipage}{0.9\textwidth}\small {\sc Figure~\thefigure. }#1\end{minipage}}

\newcommand{\textd}{\text{\rm d}\mkern0.5mu}

\newcommand{\texte}{\text{\rm e}}

\newcommand{\1}{\operatorname{\sf 1}}

\newcommand{\BB}{\mathcal B}

\newcommand{\DD}{\mathcal D}
\newcommand{\EE}{\mathcal E}

\newcommand{\HH}{\mathcal H}

\newcommand{\KK}{\mathcal K}

\newcommand{\PP}{\mathcal P}

\newcommand{\RR}{\mathcal R}

\newcommand{\BbbL}{\mathbb L}

\newcommand{\BbbP}{\mathbb P}
\newcommand{\Q}{\mathbb Q}
\newcommand{\R}{\mathbb R}
\newcommand{\BbbS}{\mathbb S}

\newcommand{\V}{\mathbb V}
\newcommand{\W}{\mathbb W}

\newcommand{\scrH}{\mathscr{H}}

\newcommand{\frG}{\fraktura G}
\newcommand{\frg}{\fraktura g}

\newcommand{\frf}{\fraktura f}
\newcommand{\frv}{{\fraktura v}_\lambda^\beta}

\def\myffrac#1#2 in #3{\raise 2.6pt\hbox{$#3 #1$}\mkern-1.5mu\raise 0.8pt\hbox{$#3/$}\mkern-1.1mu\lower 1.5pt\hbox{$#3 #2$}}
\newcommand{\ffrac}[2]{\mathchoice%
	{\myffrac{#1}{#2} in \scriptstyle}
	{\myffrac{#1}{#2} in \scriptstyle}
	{\myffrac{#1}{#2} in \scriptscriptstyle}
	{\myffrac{#1}{#2} in \scriptscriptstyle}
}

\newcommand{\ket}[1]{\left|#1\right\rangle}
\newcommand{\bra}[1]{\left\langle#1\right|}

\newcommand{\TR}{\text{\rm Tr}}

\newcommand{\betac}{\beta_c}
\newcommand{\SPAN}{\text{SPAN}}
\newcommand{\state}[1]{\left\langle #1\right\rangle}

\newcommand{\Qrob}[2]{\Q_{#1}^{#2}}
\newcommand{\Pprob}[2]{\BbbP_{#1}^{#2}}
\newcommand{\Irob}[2]{\mu_{#1}^{#2}}
\newcommand{\Wrob}[2]{\W_{#1}^{#2}}
\newcommand{\var}[2]{\V\textrm{ar}_{#1}^{#2}}
\newcommand{\lb}{\left[}
\newcommand{\rb}{\right]}

\newcommand{\lbr}{\left\{}
\newcommand{\rbr}{\right\}}

\newcommand{\la}{\left<}
\newcommand{\rab}{\right>}

\newcommand{\usigma}{\underline{\sigma}}

\newcommand{\uvarphi}{\underline{\varphi}}
\newcommand{\uxi}{\underline{\xi}}

\title[The Mean Field Quantum Ising Model]{The Phase Diagram of the Quantum Curie-Weiss Model}
\author[Chayes et. al.]{Lincoln Chayes$^1$, Nicholas Crawford$^2$, Dmitry Ioffe$^3$ \and\ Anna Levit$^3$}

\begin{document}
\thanks{\hglue-4.5mm\fontsize{9.6}{9.6}\selectfont\copyright\,2008 by L.~Chayes, N.~Crawford, D.~Ioffe A.~Levit. Reproduction, by any means, of the entire
article for non-commercial purposes is permitted without charge.\\
The research of Nicholas Crawford, Dmitry  Ioffe and Anna Levit was partly supported by the 
German-Israeli  Foundation under the grant I-870-58.6/2005.
\vspace{2mm}}

\maketitle

\vspace{-5mm}
\centerline{\textit{$^1$Department of Mathematics, University of California at Los Angeles}}
\centerline{\textit{$^2$Department of Statistics, University of California at Berkeley}}
\centerline{\textit{$^3$Department of Industrial Engineering, The Technion, Haifa, Israel}}

\begin{abstract}
This paper studies a generalization of the Curie-Weiss model (the Ising model on a complete graph) to quantum mechanics.  Using a natural probabilistic representation of this model, we give a complete picture of the phase diagram of the model in the parameters of inverse temperature and transverse field strength.  Further analysis computes the critical exponent for the decay of the order parameter in the approach to the critical curve and gives useful stability properties of a variational problem associated with the representation.
\end{abstract}

\section{Introduction}
\noindent One of the simplest classical systems exhibiting phase
transition is the Curie-Weiss model. In this model, $N$ Ising
spins
$\usigma = \lbr\sigma_i=\pm1 ;\ i=1,\dots,N\rbr$, interact via the Hamiltonian
\begin{equation}
\label{CW-hamiltonian}
H_N(\usigma)=-\frac1{2N}\sum_{i,j=1}^N\sigma_i\sigma_j,
\end{equation}
where the normalization by~$\ffrac1N$ makes $H_N$ a quantity of
order~$N$. (Generalizations to multiple-spin interactions may also
be considered; cf Sect.~\ref{Main-Results}. However, for the time
being, \eqref{CW-hamiltonian} will~suffice.)

As is well known
%~\cite{Ellis}, 
for the spins distributed according
to the measure
$\mu^{N,\beta}(\{\usigma\})\propto\texte^{-\beta
H_N(\usigma)}$, i.e., in the canonical ensemble, as~$N\to\infty$ the
law of the empirical mean $m_N(\usigma)=N^{-1}\sum_i\sigma_i$
converges to a mixture of point masses at~$\pm m_\star$ where
$m_\star=m_\star(\beta)$ is the so called spontaneous magnetization.
The phase transition in this model is manifestly seen from the
observation that $m_\star(\beta)\equiv 0$ for~$\beta\le\betac$
while~$m_\star(\beta)>0$ for~$\beta>\betac$. The function
$\beta\mapsto m_\star(\beta)$ is in fact the maximal non-negative
solution of the equation
\begin{equation}
%\label{}
m=\tanh(\beta m),
\end{equation}
while~$\betac=1$ is the so called critical temperature.

The main reason why the Curie-Weiss model is so approachable is the
fact that the Hamiltonian is, to within an additive constant, equal
to~$-\frac12 Nm_N(\usigma)^2$. This permits a very explicit
expression for the law of~$m_N$ whose concentration properties are
then readily controlled by straightforward large-deviation
arguments. However, this simple strategy breaks down once
inhomogeneous terms (i.e., those not invariant under exchanges of
the spins) are added to the Hamiltonian. One example where this
happens is the Curie-Weiss system in random external field where the
term
$\sum_i h_i \sigma_i$, with~$h_i$ sampled from an i.i.d.\ law
with zero mean, is added to~$H_N$. While rigorous analysis is still
possible in this case, the technical difficulties involved are more
substantial. Significantly more complex is the
Sherrington-Kirkpatrick version of \eqref{CW-hamiltonian}, where the
term~$\sigma_i\sigma_j$ is weighed by a (fixed) random
number~$J_{x,y}$ that has been sampled from a symmetric distribution
on~$\R$. This model possesses a beautiful underlying
structure~\cite{MPV,Parisi} which has been harnessed mathematically
only very recently~\cite{Guerra,Talagrand-book,Talagrand-paper}.

The goal of this paper is to study another natural generalization of
the Curie-Weiss model, namely to the realm of \emph{quantum
mechanics}. Here the classical spin variables~$\sigma_i$ are understood
as eigenvalues of 
$z$-component of the triplet of Pauli
matrices~$(\hat{\sigma}^{(x)},\hat{\sigma}^{(y)},\hat{\sigma}^{(z)})$---the generators
of $\mathfrak{su}(2)$---acting on the one-particle Hilbert
space~$\scrH_1=\SPAN\{\ket+,\ket-\}$. The configuration space is
replaced by the product space~$\scrH_N=\bigotimes_{i=1}^N\scrH_1$.
Classical Ising cofiguartions $\usigma\in \left\{\pm 1\right\}^N$ generate an 
orhtonormal basis $\ket{\usigma }= \otimes\ket{\sigma_i}$ of $\scrH_N$.  
If $A$ is an operator on $\scrH_1$ 
(e.g. $A = \hat{\sigma}^{(x)}$ or $A = \hat{\sigma}^{(z)}$)
then its copy $A_i$ acts on $i$-th spin component of 
%the spin operator~$\hat{\sigma}_i^{(k)}$ for the spin at~$i$ acts on a
a product vector
$\ket{\uvarphi} =\otimes_{i=1}^N\ket{\varphi_i}\in\scrH_N$ via
\begin{equation}
%\label{}
%\hat
%\sigma_i^{(k)}
A_i \ket{\uvarphi} =\ket{\varphi_1}\otimes\dots\otimes 
%\hat
%\sigma^{(k)}
A\ket{\varphi_i}\otimes\dots\otimes\ket{\varphi_N}.
\end{equation}
There are at least two natural ways to introduce quantum effects
into the Curie-Weiss model. Either one may make the interaction term
isotropic---this corresponds to the quantum Heisenberg model---or
one may consider an external transverse field. Here we focus on the latter
situation: The Hamiltonian is now an operator on~$\scrH_N$ defined
by
\begin{equation}
\label{quantum-CW-ham} \hat
\HH_N=-\frac1{2N}\sum_{i,j=1}^N\hat \sigma_i^{(z)}\hat \sigma_j^{(z)}-\sum_{i=1}^N(h\hat \sigma_i^{(z)}+\lambda\hat \sigma_i^{(x)}).
\end{equation}

%The 
Instead of 
Gibbs probability measures
 one studies 
%~$\mu^{N,\beta}$ is replaced by
 a
\emph{KMS state} $\state{-}^{N,\beta}_{\lambda ,h}$ which is a positive linear
functional on the $C^\star$-algebra of operators on~$\scrH_N$
defined by
\begin{equation}
\label{KMS}
\langle A\rangle^{N,\beta}_{\lambda ,h}=\frac{\TR( A\,\texte^{-\beta\
\HH_N})}{\TR(\texte^{-\beta \HH_N})}.
\end{equation}
As before, the parameter~$\beta$ plays the role of inverse
temperature while~$\lambda$, which corresponds to the strength of an
external field, determines the overall strength of the quantum
perturbation.

In particular, \eqref{KMS} gives rise to a probability distribution 
$\mu^{N,\beta}_{\lambda ,h}$ on classical spin configuartions
$\usigma\in\lbr\pm 1\rbr^N$, 
\begin{equation}
 \label{QGibbs}
\mu^{N,\beta}_{\lambda ,h}\lb\usigma\rb \, =\, 
\frac{\bra\usigma \texte^{-\beta \HH_N}\ket{\usigma}}
{\TR(\texte^{-\beta \HH_N})},
\end{equation}
and one can try to read a signature of phase transition in terms of concentration
properties of the latter as $N$ tends to infinity.

Despite this  relatively clean formulation, the quantum nature makes
this model very different from the classical one. Indeed, the
Hamiltonian is not simultaneously diagonalizable
with the~$\hat{\sigma}_i^{(z)}$'s for~$\lambda\ne 0$, so the very notion of ``value of
spin at~$i$'' is apparently lost. 
In fact, as we explain below, a path integral approach or, equivalently, 
a stochastic geometric representation of $ \mu^{N,\beta}_{\lambda ,h}$ 
reveals that relevant large deviation and measure concentration 
analysis in the genuine quantum case $\lambda > 0$ 
should be lifted to infinite dimensions.  
%Fortunately, a classical system may be recovered if one resorts to
%auxiliary variables.
% An approach along these lines leads to a
%graphical representation~\cite{AN} whose percolation properties are
%often related to the existence/absence of phase transitions. A study
%of a model closely related to the above---namely, its
%quantum-percolation version---based on graphical representations has
%been made recently~\cite{Ioffe}.

We refer to \cite{AKN, AN, Cam-Kl-Per, N-1, N-2} where the ideas of stochastic 
geometric representation in question were originally developed, as well as to a 
recent review \cite{ILN}. The bottom line is the representation
$\mu^{N,\beta}_{\lambda ,h}$: Let $\mathbb S_{\beta}$ be the  circle of 
circumference $\beta$.  There exists a shift invariant probability measure
$\mu^{\beta}_{\lambda }$ on piece-wise constant trajectories
 $\sigma : \mathbb S_{\beta}\mapsto \lbr\pm 1\rbr$, such that
\begin{equation}
 \label{Representation}
\mu^{N,\beta}_{\lambda ,h}\lb\usigma\rb\, \propto\, 
\bigotimes_{i=1}^{N}
\mu^{\beta}_{\lambda }
\lb
\texte^{
N\int_0^\beta [\frac12m_N(t)^2 
+h m_N (t )]
\textd t}~{\mathbf 1}_{\lbr \usigma (0) =\usigma\rbr}\rb
\end{equation}
where $\usigma (t) = (\sigma_1 (t) ,\dots ,\sigma_N (t))$ and 
$m_N (t) = N^{-1}\sum_i \sigma_i (t)$.

%Our strategy in this paper will be a variation on this idea and will
%be close in spirit to the way the classical Curie-Weiss model is
%solved. Using the Lie-Trotter formula, we will represent matrix
%elements of the Gibbs-Boltzmann weight~$\texte^{-\beta \HH_N}$ in the
%basis of classical Ising states
%$\{\ket{\sigma}\colon\sigma\in\{-1,+1\}^N\}$---the eigenbasis
%of~$\hat\sigma_i^{(z)},x=1,\dots,N$---as an expectation,
%\begin{equation}
%\label{class-repr-CW} \bra{\sigma}\texte^{-\beta
%\hat \HH_N}\ket{\tilde\sigma} =\E_{\lambda,\sigma}\bigl(\texte^{
%N\int_0^\beta [\frac12m_N(t)^2+hm_N(t)]\, \textd t}\bigr),
%\end{equation}
%over a collection $\sigma(t)=(\sigma_i(t))_{1\le x\le N}$ of~$N$
%independent Poisson point processes on circles $\mathbb S_{\beta}$
%with arrival rate~$\lambda$, one for each site $i$. The
%quantity~$m_N(t)$ is the empirical magnetization at time~$t$,
%i.e.,~$m_N(t)=N^{-1}\sum_{i=1}^N\sigma_i(t)$; the interaction is
%thus the classical interaction averaged over ``time.''  For the readers convenience, we review this %identification in the next section.

With such a  representation in place, we then apply methods of
large-deviation theory to derive the leading-order $N\to\infty$
asymptotic of these expectations.  As can be expected, many physical
properties of the quantum system may be gleaned from the properties
of the minimizer of the corresponding variational problem. 

The
variational problem (or rather dual thereof) turns out to have
intrinsic features which allow us 
bring methods of $FK$ percolation
to bear. 
 Many quantitative characteristics of the system will be
determined explicitly (at least in the limit~$N\to\infty$). In
particular, we obtain full control of the phase diagram and
the stability near optimizers for variational problem.  Moreover we
give an explicit characterization of the critical exponent for the
decay of the optimizer as one approaches the critical curve.

Let us mention existing work on issues encompassing aspects of the present paper.  The paper \cite{FSW} addresses general mean field quantum spin systems in their C$^*$ algebraic representation.  Relying on an operator theoretic version of deFinetti's Theorem (i.e. St\"{o}rmer's Theorem), they derive a `mean field equation' for the extremal states of the system and formulate variational problem which these states must solve.  More recently, the preprint \cite{Dorlas} uses a slightly different path integral representation from ours to derive a qualitative probabilistic description of the solutions to these `mean field equations'.  However, neither work attempts a detailed analysis of the Quantum Curie-Weiss Model.

The remainder of this paper is organized as follows.  In Section \ref{IntroCWMF} we review the probabilistic formulation of our problem and use this opportunity to set notation.  In Section \ref{Main-Results} we formulate the main results of this paper.  Subsequent sections are devoted to the proof of these results.
%%%%%%%%%%%%%%%%%%%%%%%%%%%%%%%%%%%%%%%%%%%%%%%%

\section{Stochastic Geometry of the Mean Field Transverse Ising Model }\label{IntroCWMF} 
Let us  explicitely describe the one-circle measure $\mu_{\beta}^{\lambda}$
in \eqref{Representation}.  
It should be noted that 
a natural approach, exploited in the papers Aizenman-Klein-Newman \cite{AKN} and Campanino-Klein-Perez
\cite{Cam-Kl-Per} among others is to realize the Gibbs state 
$\mu_{\beta}^{\lambda} = \langle\cdot\rangle_{\lambda, c}^\beta$, 
which we construct directly here, 
as the strong coupling limit of a sequence \textit{discrete} Ising systems on $\mathbb Z / n \mathbb Z$.  As a result of this approximation we get  a 
 {\em weak} 
 form of ferromagnetic 
correlation inequalities for $\langle\cdot \rangle_{\lambda, c}^\beta$ for free.  We well rely on this below without further mention.

We begin by introducing a convenient notation for expectation
values.  Given a probability measure $\mathbb P$ on a sample space
$\Omega$ and an integrable function $f: \Omega \rightarrow \mathbb
R$, let us denote the expectation value of $f$ by $\mathbb P(f)$.

Let $\mathbb P_{\lambda}^{\beta}$ be the distribution of the Poisson
point process $\xi\subset  {\mathbb S}_\beta$ of marks(which we visualize as puntucres) on the circle 
${\mathbb S}_\beta$ with
arrival intensity $\lambda$.  
Given a realization of $\xi$ let us say that a classical piece-wise
constant trajectory
 $\sigma :\BbbS_\beta\mapsto \lbr\pm 1\rbr$ is compatible with $\xi ;\
{\sigma}\sim\xi$, if  jumps of
$\sigma  (\cdot )$
 occur only at arrival times of $\xi_i$.  Note that we do not require that $\sigma$ changes
sign at each arrival of $\xi$.  Consider now a joint probability distribution, 
\begin{equation}
 \label{PhiMeasure}
\Phi^{\beta}_{\lambda ,h}\lb\textd\xi ,\textd\sigma\rb\, 
\propto
\, 
\mathbb P_{\lambda}^{\beta}\lb \textd\xi\rb \texte^{h \int_0^\beta \sigma (t)\textd t}
\delta_{\lbr \sigma\sim\xi\rbr}.
\end{equation}
Here $\delta_{\lbr \sigma\sim\xi\rbr}$ gives mass one to $\sigma$ if it is compatible with $\xi$ and gives zero mass otherwise.
In the sequel we shall supress the sub-index $h$ whenever $h=0$.  Our one-circle measure
$\mu^{\beta}_{\lambda}$ is just the $\sigma$-marginal of $\Phi_{\beta}^{\lambda }$.
Clearly, for $h\neq 0$,  $\sigma$-marginals $\mu^{\beta}_{\lambda ,h}$ of 
$\Phi^{\beta}_{\lambda ,h}$ can be recovered from $\mu^{\beta}_{\lambda}$ as 
exponential tilts by a {\em constant} magnetic field $h$, 
\[
\mu^{\beta}_{\lambda ,h}\lb \textd\sigma\rb\, \propto\, 
 \mu^{\beta}_{\lambda}\lb \texte^{h\int_0^\beta \sigma (t)\textd t}\textd\sigma\rb .
\]
Summing with respect to compatible $\sigma$-trajectories in \eqref{PhiMeasure}
we recover the  $\xi$-marginal $\Qrob{\lambda ,h}{\beta}$ of $\Phi^{\beta}_{\lambda}$:
 For each realization of $\xi$ the punctured circle $\BbbS_\beta\setminus\xi$ is a 
disjoint union of finite number  $\# (\xi)$ of connected components, 
\[  
\BbbS_\beta\setminus\xi = \bigcup_{j=1}^{\# (\xi )} I_j 
\]
Then, 
\[
 \Qrob{\lambda ,h}{\beta}\lb\textd\xi\rb\, \propto\, 
\prod_{j=1}^{\# (\xi )}\lb \texte^{h|I_j |} + \texte^{-h|I_j |}\rb 
\Pprob{\lambda}{\beta}\lb\textd\xi\rb
\]
Note that in the case $h=0$,
\begin{equation}
\label{eq:FK}
\Qrob{\lambda}{\beta}\lb \textrm{d} \xi_i\rb \, =\, \frac{2^{\#
(\xi_i )}\mathbb P_{\lambda}^{\beta}\lb \textrm{d} \xi_i\rb}{
\mathbb P_{\lambda}^{\beta}\lb 2^{\# (\xi_i )}\rb }
%= \mathbb P^{2\lambda}_{\beta}\lb  \textrm{d} \xi_i \rb 
\end{equation}
In general, 
$\lbr \Qrob{\lambda ,h}{\beta}\rbr$ is still a   family of Fortuin-Kasteleyn type 
random cluster measures.  In particular,  they possess two specific properties which we 
would like to stress:
\vskip 0.1cm

\noindent
{\em 1) Stochastic domination.} For each $h$ there exists 
an intensity
$\eta = \eta (h,\lambda )$ of the Poisson process of arrivals $\xi$ on $\BbbS_\beta$, 
such that
\begin{equation}
 \label{Domination}
\Qrob{\lambda ,h}{\beta}\, \prec\, \Pprob{\eta}{\beta}
\end{equation}
We refer to \cite{ACCN, AKN} for a startegy of proving statements 
like \eqref{Domination}, which is based on a comparison with independent 
Bernoulli percolation via the strong FKG inequality.

\noindent
{\em 2) Edwards-Sokal representation.} The spin marginal $\mu_{\lambda ,h}^{\beta}$ 
could be recovered from $\Qrob{\lambda ,h}{\beta}$ 
by the following two step procedure: First sample $\xi$
from $\Qrob{\lambda ,h}{\beta}$, and then 
independently 
paint each connected component $I$
 of
${\mathbb S}_\beta\setminus\xi$ into $\pm1$ 
with probabilities $\texte^{h|I|}/ (\texte^{h|I|}+ \texte^{-h|I|})$ and 
$\texte^{-h|I|}/ (\texte^{h|I|}+ \texte^{-h|I|})$ respectively.
\vskip 0.1cm

\noindent
In the sequel we shall repeatedly rely on the following consequence of the 
two properties above: Fluctuations of $\sigma$ under $\mu_{\lambda ,h}^{\beta}$ 
are stochastically controlled on all possible scales by arrivals of $\xi$ under 
the afore mentioned Poisson measure $\BbbP_{\eta}^\beta$.

%%%%%%%%%%%%%%%%%%%%%%%%%%%%%%%%%%%%%%%%%%%%%%%

\section{Main Results}\label{Main-Results}
A slight generalization of our Hamiltonian
of interest is given by
\begin{equation}
- \mathcal H_N= N \PP \left(\frac{1}{N} \sum_{i=1}^N
\sigma^{(z)}_i\right) + \sum_{i=1}^{N} [\lambda \sigma^{(x)}_i + 
h\sigma^{(z)}_i]  
\end{equation}
where $\PP(\cdot)$ is some polynomial 
function from 
%$\mathbb R$ 
$[-1, 1]$ to 
$\mathbb R$. 

As is apparent from \eqref{Representation}, an analysis of the phase 
diagram of the corresponding generalized CW model in
transverse field boils down to an investigation of asymptotic
properties of  the sequence of  measures
\begin{equation}
\label{eq:weighted} 
%\otimes
\Wrob{\lambda}{N, \beta}(\textrm{d}\sigma
)\, = \, \frac{\otimes\Irob{\lambda}{\beta } \lb {\exp}\lbr N
\int_0^\beta \PP(m_N (t)) \textrm{d} t\rbr\, ;\,  \textrm{d}\sigma
\rb} {\otimes\Irob{\lambda}{\beta} \lb {\exp}\lbr N \int_0^\beta
\PP(m_N (t)) \textrm{d} t \rbr \rb} ,
\end{equation}
where
\[
m_N (t)\, =\, \frac1{N}\sum_i \sigma_i (t) .
\]
This problem may be addressed with the theory of large deviations.  
\begin{theorem}
 \label{LD}
For every $h\geq 0$ and $\lambda\geq 0$ the law of 
the average $m_N(\cdot )$
 under the product measures $\otimes\Irob{\lambda ,h}{\beta }$
is exponentially tight on $L_2 (\BbbS_\beta
%, \textd t 
)$, and, furthermore it 
satisfies on the latter space a LD principle with a good rate function
$I_{\lambda ,h}^\beta $.
\end{theorem}

This follows from fairly standard methods and will be proved in the next section.
Consequently, the measures \eqref{eq:weighted} 
 are also exponentially tight on $L_2 (\BbbS_\beta 
%, \textd t 
)$ and 
they satisfy a LD principle with a good rate function
$\max_{m} \frG_{\lambda ,h}^\beta  (m) - \frG_{\lambda ,h}^\beta (\cdot )$, where, 
\[
 \frG_{\lambda ,h}^\beta  (m) \, =\, \int_0^\beta \PP(m (t))\textrm{d} t\,
-\, I_{\lambda ,h}^\beta (m).
\]
%As a result,
In particular, in  the zero magnetic field case of $h=0$,  measures \eqref{eq:weighted}  
are
 asymptotically
concentrated around solutions  of
\begin{equation}
\label{eq:var} \sup_{m}\lbr \int_0^\beta \PP(m (t))\textrm{d} t\,
-\, I_{\lambda }^\beta (m)\rbr\, \equiv  \, \sup_{m} \frG_{\lambda }^\beta  (m) .
\end{equation}
%where $I$ is the large deviation rate function for the average $m_N$
 %under the product measures $\otimes\Irob{\lambda}{\beta }$. 
Below we shall use $(\cdot ,\cdot )_\beta$ to denote the scalar product
in $L_2 ({\mathbb S}_\beta )$ and $\|\cdot \|_\beta$ to denote the corresponding
norm.  Then, 
since  we formulate
 our large deviation principle in $L_2 ({\mathbb S}_\beta )$, 
%, \textrm{d} t 
%)$, then,
%using $(\cdot ,\cdot )_\beta$ to denote the corresponding scalar product,
\begin{equation}
\label{eq:ILambda} I^{\beta}_{\lambda ,h} (m)\, =\, \sup_{g}\lbr (g, m)_\beta - 
\Lambda_{\lambda ,h}^\beta 
(g)\rbr
\end{equation}
where
\begin{equation}
\Lambda _{\lambda ,h}^\beta (g)\, =\,
\log\Irob{\lambda ,h}{\beta}\lb {e}^{(g,\sigma)_\beta}\rb.
\end{equation}
For the solution of the variational problem and, in the computation of the corresponding transitional curve, one must proceed on a case by case basis.  Let us specialize to the \textit{quadratic} problem.

In order to study phase coexistence in this model, we must consider the case
of zero magnetic field in $z$-direction, $h=0$. 
Let us define
\begin{equation}
\label{eq:CWtcurve} \frf (\lambda ,\beta)\, = \,
\frac1{\beta}\var{\lambda}{\beta}\lb (\sigma ,\mathbf 1 )_\beta\rb
\end{equation}
where $\var{\lambda}{\beta}$ is the variance under
 the one-circle spin measure $\Irob{\lambda}{\beta}$.
A straightforward calculation (which we present in
 Section \ref{sec:1D}) shows that
\begin{equation}
\frf(\lambda, \beta)
=\, \frac1{\lambda}\tanh (\lambda\beta ).
\end{equation}

For the statement of our main theorem, 
let 
\begin{equation}
s_4^h\left(\lambda, \beta \right) = - \frac{\textd^4}{ \textd h^4}
%{\Big\vert}_{h=0} 
\Lambda_\lambda^\beta \left(h \cdot \bf 1\right)\quad\text{and}\quad 
s_4 \left(\lambda, \beta \right)= s_4^0\left(\lambda, \beta \right) .
\end{equation}
Our analysis below will show that there exists $C=C(\lambda , \beta)>0$ so that
\begin{equation}
 \frac{\textd^3}{ \textd h^3}{\Big\vert}_{h=h_0} 
\Lambda_\lambda^\beta \left(h \cdot \bf 1\right) \leq - C h_0.
\end{equation}
for $h_0$ close to zero.
Thus, in view of the spin-flip symmetry and continuity of 
semi-invariants $s_4^h\left(\lambda, \beta \right)$ 
 in $h$,   $s_4\left(\lambda, \beta \right)$ is strictly positive.

The results of our analysis are summarized as follows:
\begin{theorem}\label{Main}
Let $\PP(x)= \frac{1}{2} x^2$.  The variational problem
\eqref{eq:var} at zero magnetic field $h=0$ 
has constant maximizers   $\pm m^* (\lambda ,\beta
)\!\cdot\! \mathbf 1$, where the spontaneous $\textrm{z}$-magnetization
$m^*$ satisfies:
\begin{enumerate}
\item If $\frf (\lambda ,\beta )\leq 1$, then $m^* = 0$.
\item If $\frf (\lambda ,\beta ) > 1$, then $m^* >0$, and, consequently
there are two distinct solutions to \eqref{eq:var}.
\end{enumerate}
Furthermore, away from the critical curve the solutions $\pm
m^*\!\!\cdot \!\mathbf 1$ are stable in the following sense: There exists
$ c = c (\lambda ,\beta )>0$
 and a strictly convex symmetric function $U$ with 
 a $U (r)\sim r\log r$ growth at infinity, 
%$U (0) = 0$ and
%$U^{\prime\prime} (0) >0$
such that
\begin{equation}
\label{eq:stable} 
\frG_\lambda^\beta  (\pm m^*\!\!\cdot\! \mathbf 1 ) - \frG_\lambda^\beta (m)\, \geq\,
c\max \lbr \| m \pm m^*\!\!\cdot\! \mathbf 1\|_{
%2;\;
\beta}^2\,  , \int_0^\beta U(m^\prime (t) )\textrm{d} t \rbr ,
\end{equation}
where, 
\[
 \| m \pm m^*\!\! \cdot \!\mathbf 1\|_{
%2;\;
\beta}^2 \, =\, 
\min\lbr \| m - m^*\!\! \cdot\!\mathbf 1\|_{
%2; \;
\beta}^2\,  ,\,   \| m +m^*\!\!\cdot \!\mathbf 1\|_{
%2; \;
\beta}^2\rbr. 
\]
As a result, the variational 
problem \eqref{eq:var} is  also stable in the supremum norm $\|\cdot\|_{\rm sup}$.
Namely, 
 there exists a constant $c_{\rm sup}  = c_{\rm sup}(\lambda ,\beta ) < \infty $, such that
\begin{equation}
 \label{supbound}
\frG_\lambda^\beta  (\pm m^*\!\!\cdot\! \mathbf 1 ) - \frG_\lambda^\beta (m)\, \geq\, 
{\rm exp}\lbr -\frac{c_{\rm sup}}{ \|m \pm m^*\!\! \cdot \!\mathbf 1\|_{
%2;\;
\rm sup}}\rbr .
\end{equation}

\noindent
Finally we have the following expression 
for the decay of $m^* >0$ in the super-critical region near
the critical curve:
\begin{equation}
\label{eq:decay}
m^* = m^* (\lambda ,\beta ) = 
\sqrt{\frac{ 6 \beta \left(\frf (\lambda ,\beta )  - 1\right)}{s_4 (\lambda, \beta)}}\left(1 + o(1)\right)
\end{equation}
% \biggl \lvert 1-\frac{1}{\beta}\textrm{$\mathbb
%V$ar}_0\bigl[ (\sigma, \mathbf 1)_\beta\bigr] \biggr \rvert
where the implicit constants depend on $\beta$ and $\lambda$ but
are bounded below in compact regions of the parameter space.
\end{theorem}

\begin{remark}
We remark that the second term of \eqref{eq:stable} is particularly 
important in
the super-critical regime ($\frf (\lambda ,\beta ) > 1$) since it
rules out trajectories of $m_N (\cdot )$ with rapid transitions
between the optimal values $\pm m^*$.

\end{remark}

\begin{remark}
We remark that there are alternate approaches to quantum Ising systems in which the random cluster representation (and percolation therein) play a more central role, e.g.~as in \cite{Ioffe} \cite{ILN} and  \cite{Grim}.  While these ideas will provide a backdrop for our derivation, at the core we shall adhere to a more traditional thermodynamic/large deviations approach:  Various results, e.g~stability, are stronger in this context.  As a bit of foreshadowing to future work, we remark that the spin system formulation is at present essential for the extension to finite dimensional systems.
On the other hand we would like to stress that our results automatically imply that, as it 
was claimed  in \cite{Ioffe} and later on conjectured in \cite{Grim}; 
 $\frf (\lambda ,\beta ) =1$,  is the equation of the critical curve for the $q=2$ quantum 
 FK-percolation on complete graphs.
\end{remark}
In Section~\ref{LDProduct} we shall prove Theorem~\ref{LD}. 
Then the remainder of this paper
develops our proof of Theorem~\ref{Main}.  

\begin{comment}
In the concluding Section~\ref{FK}
we discuss implications for the quantum FK model on complete graphs.
\end{comment}
%%%%%%%%%%%%%%%%%%%%%%%%%%%%%%%%%%%%%%%%%%%%%%%%%
\section{Large deviations under product mesures}
\label{LDProduct}
In this section we study asymptotic properties of the average $m_N (\cdot )$ 
under the sequence of product mesures $\otimes\mu_{\lambda ,h}^{\beta}$.  Let us begin with the proof of Theorem \ref{LD}.

\begin{proof}[Proof of Theorem \ref{LD}]
Let us start with the exponential tightness on $L_2 (\BbbS_\beta )$.  Recall
the following realization of compact subsets of $L_2 (\BbbS_\beta )$:
 Let $a :[0,\beta/2 ]\mapsto \R_+$ be a continuous non-decreasing function with $a(0)=0$.
Then, 
\begin{equation}
 \label{Kset}
\KK_a\, =\, \lbr g\in L_2 (\BbbS_\beta )~:~\|g\|_{
%2;\; 
\beta}\leq 1\, \text{and}\, 
\| g(s +\cdot )-g (\cdot )\|_{
%2; \;
 \beta} \leq a(s) ~\forall~s\in [0,\beta/2 ]\rbr
\end{equation}
is compact. Now, in view of the Edwards-Sokal representation, 
\[
\| m_N(s +\cdot )-m_N (\cdot )\|_{
%2; \; 
\beta} \, \leq\, 2s\frac1{N}\sum_{i=1}^N\xi_i (\BbbS_\beta ) .
\]
The latter expression is increasing in $\uxi =(\xi_1, \dots ,\xi_N )$, hence it is 
exponentially bounded
in view of   stochastic domination by product Poisson measures 
$\otimes\BbbP_{\eta}^{\beta}$. %\qed

\noindent
Large deviations with good rate function $I_{\lambda ,h}$ defined in 
\eqref{eq:ILambda} is a standard consequence \cite{Baldi}. % \qed
\end{proof}

\noindent
In the sequel we shall need the following stability estimate which implies that
product mesures $\otimes\mu_{\lambda ,h}^{\beta}$ are sharply concentrated
around constant functions.
\begin{theorem}
 \label{thm:IStability}
Let $\lambda ,h\geq 0$ be fixed, and let $\eta  $ be an intensity of 
the Poisson process of arrivals on $\BbbS_\beta$ such that \eqref{Domination} holds.
Then, 
\begin{equation}
\label{IStability}
 I^\beta_{\lambda ,h} (m)\, \geq\, 
\begin{cases}
 \infty ,\, &\text{if $m$ is not absolutely continuous}\\
\int_0^\beta U_\eta (m^\prime (t))\textd t ,\, &\text{otherwise} ,
\end{cases}
\end{equation}
 where $U_\eta$ 
is an 
even smooth strictly convex function with
% a quadratic minimum at zero 
%and 
 a superlinear growth at infinity. 
\end{theorem}
\begin{remark}
 The function $U_\eta$ is defined in \eqref{Ueta} below. Note that
it  is negative in a neighbourhood of the origin.  What is important, however, 
is a super-linear growth at infinity which, as we shall see in the 
concluding Section~\ref{sec:FullStability},  controls rapid oscillations for the 
original variational problem \eqref{eq:var}.
\end{remark}

\begin{proof}
Let $\RR= \lbr 0 < t_1 <t_2 <\dots <t_n=t_0 <\beta\rbr$ be a 
partition of $\BbbS_\beta$.  Consider 
an $n$-dimensional random vector
\[
 z_N^{\RR}\, =\, \lb m_N (t_1 )- m_N (t_0), m_N (t_2 ) -m_N (t_1), 
\dots ,m_N (t_n )- m_N (t_{n-1})\rb .
\]
Evidently, $z_N^{\RR}$ satisfies a LD principle with rate function $I^{\RR}$, 
\begin{equation}
 \label{IRR}
I^{\RR} (z_1 ,\dots ,z_n )\, =\, \max_{g_1 ,\dots ,g_n}\lbr
\sum_i g_i z_i - \Lambda^{\RR} (g_1 ,\dots ,g_n )\rbr ,
\end{equation}
where
\begin{equation}
 \label{LRR}
\Lambda^{\RR} (g_1 ,\dots ,g_n )\, =\, \log \mu_{\lambda ,h}^{\beta}
\lb \texte^{\sum_i g_i (\sigma (t_{i}) -\sigma ( t_{i-1}))}\rb 
\end{equation}
By contraction \footnote{This is a somewhat formal argument which, however, 
could be easily turned into a rigorous one using a standard molifying procedure}, 
\begin{equation}
 \label{contraction}
I^\beta_{\lambda , h} (m)\, \geq\, I^{\RR} \lb m(t_1) - m(t_0 ), \dots ,m(t_n )- m(t_{n-1}) \rb .
\end{equation}
It, therefore, remains to derive an appropriate lower bound on $I^{\RR}$. 
Alternatively, we may seek for an upper bound on $\Lambda^{\RR}$.  Notice that 
in view of the Edwards-Sokal representation, 
\[
 \texte^{\sum_i g_i (\sigma (t_{i}) -\sigma ( t_{i-1}))}\, \leq\, 
\prod_i \lb 1 + \mathbf{1}_{\lbr \xi ([t_{i-1}, t_i ]) >0\rbr}(\texte^{2g_i}+\texte^{-2g_i})\rb .
\]
At this stage we are delighted since the latter expression is monotone in $\xi$.  Hence, 
by virtue of \eqref{Domination}, 
\[
 \mu_{\lambda ,h}^{\beta}\lb 
 \texte^{\sum_i g_i (\sigma (t_{i}) -\sigma ( t_{i-1}))}\rb \, \leq\, 
\prod_i\lb 1+ (1 - \texte^{\eta (t_{i-1} - t_i)}) (\texte^{2g_i}+\texte^{-2g_i}) \rb .
\]
Consequently, 
\[
 \Lambda^{\RR} (g_1 ,\dots ,g_n )\, \leq\, \eta \sum_i |t_{i}-t_{i-1}| H (g_i ) ,
\]
where $H (g) = (\texte^{2g } + \texte^{-2g })$.  By duality, 
\[
 I^{\RR} (z_1 ,\dots ,z_n )\, \geq\, \sum_i |t_{i}-t_{i-1}| U_\eta\lb \frac{z_i}{|t_{i}-t_{i-1}|}\rb ,
\]
with, 
\begin{equation}
 \label{Ueta}
U_\eta (z)\, =\, \eta H^* \lb \frac{z}{\eta }\rb ,
\end{equation}
with $H^*$ being the Legendre transform  of $H$. Thus, $U_\eta$ is an 
even smooth strictly convex function.
%with a quadratic minimum at zero. 
Furthermore $U_\eta (m) \sim  |m|\log |m|$ at infinity. Note, however, that $U_\eta$ is 
negative in a neighbourhood of the origin, in particular $U_\eta (0) = -2\eta$.  Going back to 
\eqref{contraction} we obtain that the rate function $I_{\lambda , h}^{\beta}$
satisfies
\[
 I_{\lambda , h}^{\beta} (m) \, \geq\, \sum_i |t_{i}-t_{i-1}| U_\eta\lb 
\frac{m(t_i)-m(t_{i-1})}{|t_{i}-t_{i-1}|}\rb ,
\]
 for every partition $\RR$ of $\BbbS_\beta$.  Taking supremum over $\RR$-s we 
recover \eqref{IStability}
\end{proof}
\newpage

\section{Duality and Reduction to One Dimensional Variational Problem}
\label{S:R-P}
Recall that we are restricting attention to the case of quadratic interaction
$\PP (x) = x^2 /2$. 
In particular, the functional
\[
 \Psi (m)\, =\, \frac12\int_0^\beta m^2 (t)\textd t\, =\, \frac12\| m\|^2_{
%2; \;
\beta} 
\]
is convex and Gateaux differentiable on $L_2 (\BbbS_\beta,\textrm{d} t )$.
Since, as we already know, $I_{\lambda}$ has compact level sets, 
the supremum in \eqref{eq:var} is attained.  Let $\bar{m}$ be a 
solution,
\[
 %\frac12
\Psi (\bar m) - I^\beta_{\lambda
%, h
} (\bar m)\, =\, \max_m\lbr  
%\frac12
\Psi ( m) - I^\beta_{\lambda
%, h
} ( m)\rbr
\]
Then, 
\[
 I^\beta_{\lambda
%, h
} (m ) - I^\beta_{\lambda
%, h
} (\bar{m})\, \geq\, \Psi (m) - \Psi (\bar{m})\, \geq\, 
( m -\bar{m} ,\bar{m})_\beta , 
\]
for every $m\in L_2 (\BbbS_\beta
%, \textrm{d} t
)$. 
It follows that $I_{\lambda}^\beta $ is sub-differential
at $\bar{m}$ and, moreover, $\bar{m}\in\partial I^\beta_{\lambda
%, h
} (\bar{m})$.
But that means 
$\Lambda_\lambda^\beta  (\bar{m})  + I^\beta_{\lambda
%, h
} (\bar{m}) = \|\bar{m}\|^2_{
%2;
\beta}$.
Hence,  for every $m$ and $h$, 
\[
 \Lambda_\lambda^\beta (\bar{m})  - \frac12\|\bar{m}\|_\beta^2\, \geq\, 
\Psi (m ) - I^\beta_{\lambda
%, h
} (m) = \lbr (h,m)_\beta -I^\beta_{\lambda
%, h
} (m) \rbr - 
\lbr (h,m)_\beta - \Psi (m) \rbr .
\]
We conclude: If $\bar{m}$ is a solution of \eqref{eq:var}, then
\begin{equation}
 \label{eq:dualvar}
\Lambda_\lambda^\beta (\bar{m})  - \frac12\|\bar{m}\|^2_{
%2;
\beta
}\, =\, \max_h\lbr 
\Lambda_\lambda^\beta  ( h)  - \frac12\| h\|^2_{
%2;
\beta}\rbr 
\end{equation}
Since $\Lambda_\lambda^\beta $ is Gateaux differentiable and, in particular, it is 
everywhere sub-differentiable the reverse conclusion is equally true: if $\bar{m}$ is 
a solution to \eqref{eq:dualvar} then it is a solution of 
\eqref{eq:var}.  Note that far reaching generalizations of the above arguments may be 
found in \cite{Toland}.   An application for general interactions $\PP$ will be attended
elsewhere.  Here we continue to stick to the quadratic case. Also, since some of the computations
below will be done in a greater generality, 
we shall casually rely on the fact that $\bar{m}$ is a solution to 
\eqref{eq:var} if and only if it is a solution of the dual problem \eqref{eq:dualvar}.
The latter  happens to be more susceptible to analysis.
\vskip 0.1cm 

\noindent
Our next step is to assemble facts
which reduce \eqref{eq:dualvar}  to a study of constant fields $h = c \cdot
\mathbf 1$ for $c \in \mathbb R$.  In what follows 
%does not depend on
%the 
%whether 
we rely only on the fact that the polynomial function  $\PP$
% is quadratic. 
has
 a
super-linear  growth at infinity. In particular, classical Legendre 
transform $\PP^*$ is a well defined function on $\R$. 
 Let integral functional $\Psi: L_2 (\BbbS_\beta
%, \textd t 
) \rightarrow
\mathbb R$ be defined by
\begin{equation}
\Psi(m) = \int_0^\beta \PP(m_t) \textd t.
\end{equation}
and for $h \in L^2(\mathbb S_\beta, \textd t)$ let 
\begin{equation}
\label{Psistar}
\Psi^*(h)= \sup_{m \in \BB} \left\{ (h,m)_\beta - \Psi(m) \right\}\, =\, \int_0^\beta \PP^*(h_t) \textd t
\end{equation}
 be its Fenchel transform.

\begin{lemma}\label{min-cons}
For any polynomial $\PP(x)$ of a super-linear growth at infinity 
, let $\Psi$ and $\Psi^*$ be defined as
above.
%  Let $\Lambda(h)$ denote the \footnote{name changed to log-moment} 
%log-moment generating function of
%the free measure $\textd \mathbb Q_{\lambda, \beta}$. 
%the one-circle measure $\mu_\lambda^\beta$\footnote{correction for measure notation}
 Then we have 
\begin{equation}
\label{ConstBound}
\sup_{h \in L_2 (\mathbb S_\beta
%, \textd t
)}\lbr  \Lambda_\lambda^\beta  (h)
 - \Psi^*(h)\rbr = \sup_{c\in
\mathbb R} \lbr \Lambda_\lambda^\beta  (c \cdot \mathbf 1) - \beta \PP^*(c)\rbr .
\end{equation}
\end{lemma}
In particular this lemma implies that the constant
function $c\cdot \mathbf 1$ optimizes the left hand side if and
only if $c$ optimizes the right hand side.

To prove Lemma \ref{min-cons}, we need 
 the following 
%to make a few 
observation: 
For any $h \in L_2 (\mathbb S_\beta)$, let $h_l$ and $h_r$ be the pair
functions obtained by reflecting $h$  about $0$ and $\frac{\beta}{2}$. That
is let
\begin{equation}
h_l(t)=
\begin{cases}
h(t) & \text{ if $t \in [0, \frac{\beta}{2}]$} \\
h(\beta - t) & \text{ if $t \in [\frac{\beta}{2}, \beta]$}
\end{cases}
\end{equation}
and $h_r$ be defined analogously.
% Let $\mathcal S$ denote the
%collection of functions in $L^2(\mathbb S_\beta)$ which are symmetric
%about $\ffrac \beta2$.
\begin{lemma}\label{RP-for-fields}
Let $h \in L_2(\mathbb S_\beta)$ and $h_l, h_r$ be as above. Then we have
\begin{equation}\label{subadd}
\Lambda_\lambda^\beta  (h) \leq \frac{1}{2} \left[\Lambda_\lambda^\beta   (h_l) +
\Lambda_\lambda^\beta   (h_r)\right].
\end{equation}
Consequently,  for every $h\in L_2 (\mathbb S_\beta)$, 
\begin{equation}
 \label{IntInequality}
\Lambda_\lambda^\beta   (h)\, \leq\, \frac{1}{\beta}\int_0^\beta 
\Lambda_\lambda^\beta   
(h(t)\!\cdot\!  \mathbf 1) \textd t .
\end{equation}
\end{lemma}
We note that Lemma \ref{RP-for-fields} was originally proved in a
somewhat more general (but completely analogous) context in
\cite{Dorlas}.
\begin{proof}
Inequality \eqref{subadd} follows from reflection positivity of classical 
ferromagnetic Ising systems: 
Recall that 
%the measure space 
$
%\left(\Omega, 
%\textd \mathbb
%Q_{\lambda, 
\mu_\lambda^\beta
%\right)
$ is the continuum limit of a sequence of
discrete Ising models with asymptotically singular interactions. Let
$M \in \mathbb N$ be fixed and suppose that the field $h\in \DD$, that  is
 $h$ is 
piecewise constant on the dyadic intervals $[\frac{k}{2^M} \beta,
\frac{k+1}{2^M} \beta)$. 
Expressing the discrete Ising model on a
circle via transfer matrices and using the generalized Cauchy-Schwartz
inequality \cite{FILS1}, one may check equation \eqref{subadd} directly whenever
the field $h\in \mathcal D$.  Passing to limits with respect to the
discrete Ising models, we immediately have \eqref{subadd} for any
dyadic piecewise constant function $h$.  Standard approximation
arguments extend the bound to the general case.  Finally, 
\eqref{IntInequality} follows from repeated applications of \eqref{subadd} 
\end{proof}
%We remark that, at bottom, the previous proof relies on the generalized H\"{o}lder inequality from \cite{FILS1}, \footnote{added}that is on reflection positivity.
Let us go back to the proof of Lemma \ref{min-cons}.  By
 \eqref{Psistar} and \eqref{IntInequality}, 
\begin{equation}
 \label{Dorlas}
 \Lambda_\lambda^\beta   (h ) - \Psi^* (h)\, \leq\, \int_0^\beta \left( \frac{1}{\beta}
\Lambda_\lambda^\beta    (h (t)\!\cdot\! \mathbf 1) - \PP^* (h (t))\right)
\textd t  
\end{equation}
\eqref{ConstBound} follows.\qed
\newpage

\section{Analysis and stability of the one-dimensional problem}
\label{sec:1D}
Following \eqref{ConstBound}  we need to analyze the following 
one-dimensional problem, 
\begin{equation}
 \label{eq:1D}
\sup_c \frg_\lambda^\beta (c)\, \equiv\, 
\sup_{c\in \R}\lbr \Lambda_\lambda^\beta  (c\!\cdot\! \mathbf{1} ) - \frac{\beta c^2}{2}\rbr .
\end{equation}
%Let us assume for this section that $\PP(x)= \frac{x^2}{2}$.  In
%this case.
Evidently, $|\Lambda_\lambda^\beta (c\!\cdot\! \mathbf{1} )|\leq\beta |c|$. Hence, 
\begin{equation}\label{growth}
\lim_{|c |\rightarrow \infty}\Lambda_\lambda^\beta (c\!\cdot\! \mathbf 1) - \frac{\beta c^2}{2} = - \infty
\end{equation}
so that solutions to \eqref{eq:1D} always exist.  
Any such solution $c$ is a critical point satisfying
%\begin{proof}[Proof of Theorem \ref{Main}:  The Critical Curve]
%We look for nontrivial solutions to the equation
\begin{equation}\label{const-MFE}
c \, = \, \frac{1}{\beta} \frac{\textd }{\textd c}\Lambda_\lambda^\beta  (c\!\cdot\!\mathbf{1} )\, 
=\, \frac{1}{\beta} 
\mu_{\lambda , c}^{\beta}\lb (\sigma ,\mathbf{1})_\beta\rb\, \equiv\, M (c) 
\end{equation}
By the FKG and GHS inequalities (cf. the discussion at the beginning of Section $2$) the function $\frac{\textd}{\textd
c} M(c)
%\big{|}_{h=u} 
= \frac{1}{\beta}\textrm{$\mathbb V$ar}^\beta_{\lambda, c} \left((\sigma,
\mathbf 1)_\beta  \right) $ is positive and decreasing in $c$. 
% One
%easily checks that it is continuous using the bounded convergence
%theorem. 
In fact, as we shall check below it is strictly decreasing on $[0,\infty )$. 
Therefore Equation \eqref{const-MFE} has nontrivial
solutions if and only if
\begin{equation}
\label{TheCondition}
\frac{\textd}{\textd c} M(c)\big{|}_{c=0} = \frac{1}{\beta}
\textrm{$\mathbb{V}$ar}^\beta_{\lambda}\left[\left(\sigma, \mathbf
1 \right)_\beta\right] \, >\, 1 .
\end{equation}
%where the variance is taken with respect to the free measure
%$\frac{\textd \mathbb Q^{\lambda}_{\beta}(\sigma)}{Z(0)}$. 
 As $\bigl \langle
\left(\sigma, \mathbf 1 \right)_\beta\bigr \rangle_{\lambda}^{\beta}= 0$ by spin
flip symmetry, $\textrm{$\mathbb
V$ar}^\beta_{\lambda}\left[\left(\sigma, \mathbf 1
\right)_\beta\right] = \mu^\beta_{\lambda }\left( \left(\sigma, \mathbf 1
\right)_\beta^2 \right)$.

We are reduced to computing
\begin{equation}
\mu^\beta_{\lambda}\left(\left(\sigma, \mathbf 1 \right)_\beta^2\right)=
\beta \int_0^\beta \mu^\beta_{\lambda}\left(\sigma_0 \sigma_t\right)
\textd t,
\end{equation}
where we have used rotational invariance in the equality. In view of the Edwards-Sokal 
representation, the
$ \pm 1$ symmetry implies
\begin{equation}
 \mu^\beta_{\lambda}\left(\sigma_0 \sigma_t\right)
 = \frac{\mathbb
P_\lambda^\beta\lb  0 \leftrightarrow t~;~2^{\#(\xi )}\rb }{\mathbb
P_\lambda^\beta\left[2^{\# (\xi )}\right]}
\end{equation}
where, as in \eqref{eq:FK} before,  $\#(\xi )$ denotes the number of connected components
determined by the underlying Poisson process and $\left\{0\leftrightarrow
t\right\}$ denotes the event the $0$ and $t$ are in the same connected
component of the complement of the arrival points.

There are two cases to consider: either there is an arrival in $[t,
\beta]$ or there is not. Taking this under consideration, we have
\begin{equation}\label{conn}
\mathbb P_\lambda^\beta \lb  0 \leftrightarrow t ~;~2^{\#(\xi )}\rb = e^{\lambda (\beta
-2t)}+ e^{-\lambda (\beta-2t)}.
\end{equation}
A simple computation shows
\begin{equation}
\mathbb P_\lambda^\beta\left[2^{\# (\xi )}\right] = e^{-\lambda
\beta} + e^{\lambda \beta}
\end{equation}
Upon integrating, we find
\begin{equation}
\frac{1}{\beta} \textrm{$\mathbb
V$ar}^\beta_{\lambda}\left(\left(\sigma, \mathbf 1
\right)_\beta\right) = \frac{1}{\lambda} \tanh(\lambda \beta).
\end{equation}
%\end{proof}
Conclusions (1) and (2) of Theorem~\ref{Main} follow as soon as we show that 
\[
\frv (c)\, \equiv\, 
%\frac{1}{\beta}
\frac{1}{\beta} \textrm{$\mathbb
V$ar}^\beta_{\lambda, c}\left(\left(\sigma, \mathbf 1
\right)_\beta\right)
\]
is indeed strictly decreasing on $[0,\infty )$. To this end,  note that, 
\begin{equation}
\label{Eq:3rd-Der}
 \frac{\textd \frv (c)}{\textd c}\, =\, \int_0^\beta\int_0^\beta\int_0^\beta
U_{\lambda, c}^\beta (r,s,t)\textd r\textd s\textd t ,
\end{equation}
where, 
 \begin{multline}
U_{\lambda, c}^\beta (r,s,t) =  \bigl \langle \sigma(r) \sigma(s) \sigma(t)\bigr
\rangle_{\lambda, c}^\beta-  \bigl \langle \sigma(r) \sigma(s) \bigr 
\rangle_{\lambda, c}^\beta
\bigl
\langle \sigma(t) \bigr \rangle_{\lambda, c}^\beta
- \bigl \langle
\sigma(s) \sigma(t) \bigr \rangle_{\lambda, c}^\beta 
\bigl \langle \sigma(r) \bigr \rangle_{\lambda, c}^\beta \\
 -\bigl \langle \sigma(r) \sigma(t) \bigr \rangle_{\lambda, c}^\beta \bigl \langle \sigma(s) 
\bigr \rangle_{\lambda, c}^\beta
+ 2 \bigl \langle \sigma(r)\bigr \rangle_{\lambda, c}^\beta \bigl \langle \sigma(s)
\bigr \rangle_{\lambda, c}^\beta \bigl \langle \sigma(t) \bigr \rangle_{\lambda, c}^\beta
\end{multline}
the third Ursell function for $\mu_{\lambda, c}^\beta$. 
 %the circle measure weighted by the
%external field $h$.  Note that, in a formal sense, this is the third derivative of
%$\Lambda(h)$).
%A natural approach, exploited in the papers Aizenman-Klein-Newman \cite{AKN} and Campanino-Klein-Perez
%\cite{Cam-Kl-Per} among others is to realize the Gibbs state $\langle\cdot\rangle_{\lambda, c}^\beta$ as the strong coupling limit of a sequence \textit{discrete} Ising systems on $\mathbb Z / n \mathbb Z$.  As a result of this approximation we get 
 %{\em weak} 
%correlation inequalities for $\langle\cdot \rangle_{\lambda, c}^\beta$ for free. 
In view of the above mentioned identification of $\mu_{\lambda ,c}^{\beta}$
as a limit of one-dimensional ferromagnetic Ising models, 
 we state the following result without proof: 
%\begin{lemma}\label{FKG-GHS}[FKG and GHS inequalities]
For any 
$c\geq 0$, 
%positive, bounded measurable function $h$ on $\Omega$, the
%weighted measure
%\%begin{equation}
%\frac{\textd \mathbb Q^{\lambda}_\beta(\sigma)}{Z(h)}
%\end{equation}
%has the FKG property with respect to the partial order defined
%above. Moreover, we have
\begin{equation}
U^\beta_{\lambda, c}(r,s,t) \leq 0 \text{\: for all $r,s,t \in \mathbb S_\beta.$}
\end{equation}
%\end{lemma}
In order to derive a strict bound we shall employ a 
similar approximation, which will allow us to transfer the random current
representation used by Aizenman \cite{Aizenman} and
Aizenman-Fernandez \cite{AF} to this continuous time setting.
We obtain the following lemma as a result.  
%This will be needed for the stability estimate
%\eqref{eq:stable} and the identification of the critical exponent. 
 While the proof is not difficult, it is tedious and so we relegate it to the appendix.
\begin{lemma}\label{RCB}[Random Current Bound]
For any positive constant field $c$ and any triplet of points $r,s,t
\in \mathbb S_\beta$, we have the bound
\begin{equation}\label{thebound}
U_{\lambda, c}^\beta(r,s,t) \leq - \bigl \langle \sigma_0 \bigr \rangle_{\lambda, c}^\beta e^{-(4
\lambda+ 2c) \beta} \frac{\lambda^2}{2}\left[ (s-r)^2 +
(t-s)^2 + (\beta+r-t)^2\right].
\end{equation}
\end{lemma}
%\eqref{thebound} readily implies that either \eqref{const-MFE} has only trivial 
%solution $c=0$, or it has three solutions
\begin{remark}
We note further that one can show that
\begin{equation}\label{thebound-1}
U_{\lambda, c}^\beta(r,s,t) \geq \\ -C(\lambda, \beta)\left(\lambda \beta\right)^2 \bigl \langle \sigma_0 \bigr \rangle_{\lambda, c}.
\end{equation}
for some $C>0$, bounded away from zero on compact regions of parameter space.
This will not play any role below and so we do not pursue the bound here.
\end{remark}

Thereby, the phase diagrame \eqref{TheCondition} is justified. 
Let us turn to the stability issue: As  follows from
\eqref{Eq:3rd-Der} and \eqref{thebound},  there exists $\chi = \chi (\lambda ,\beta )$, such that for any $c\geq 0$, 
\[
 \frac{\textd \frv (c)}{\textd c}\, \leq\, - \chi  \langle \sigma_0 \bigr \rangle_{\lambda, c}^\beta .
\]
As a result, 
for any $0\leq c_1 < c_2$, 
\begin{equation}
 \label{vbound}
\chi \langle \sigma_0 \bigr \rangle_{\lambda, c_1}^\beta (c_2 -c_1) \, \leq\,  
\frv (c_1 ) - \frv (c_2 )\, \leq\, 
\chi \langle \sigma_0 \bigr \rangle_{\lambda, c_2}^\beta (c_2 -c_1) 
\end{equation}
Assume that \eqref{TheCondition} holds, and let $c^* = m^* >0$ be the positive maximizer of $\frg_\lambda^\beta$.  Then, of course,  $0 = \frac{\textd}{\textd c}\frg_\lambda^\beta  (c^*)$. Let us derive an 
upper bound on 
\[
\frac{\textd}{\textd c}\frg_\lambda^\beta (c)\, =\, \beta c\lb \frac1{c}\int_0^c\frv (\tau )\textd \tau - 1\rb.
\]
For $c > c^*$, 
\[
 \frac1{c}\int_0^c\frv (\tau )\textd \tau\, =\, 
\frac1{c^*}\int_0^{c^*}\frv \lb \tau \frac{c}{c^*}\rb \textd \tau\, =\, 
1 + \frac1{c^*}\int_0^{c^*}\lbr \frv \lb \tau \frac{c}{c^*}\rb - 
\frv (\tau )\rbr \textd \tau .
\]
By  \eqref{vbound} the second term  is bounded above by
\begin{multline}
\label{cgeqcstar}
 -\, \chi \int_0^{c^*} \langle \sigma_0 \bigr \rangle_{\lambda, \tau}^\beta\lb\frac{c}{c^*} -1\rb
\tau\textd\tau\,\\ =\,  -\chi (c -c^* )\frac1{c^*}
\int_0^{c^*} \langle \sigma_0 \bigr \rangle_{\lambda, \tau}^\beta\tau\textd\tau\, 
\, \leq\, 
-\frac{\chi (c^*)^2}{3} (c -c^* ) ,
\end{multline}
where at the last step we have relied on $\langle \sigma_0 \bigr \rangle_{\lambda, \tau}
 \geq \tau$ for every
$\tau\leq c^*$. 

\noindent
In a completely similar fashion we deduce that, 
\begin{equation}
 \label{cleqcstar}
\frac1{c}\int_0^c\frv (\tau )\textd \tau - 1\, \geq\, \frac{\chi (c^*)^2}{3} (c^* - c ) ,
\end{equation}
whenever, $c \leq c^*$.  We have proved:
\begin{lemma} Assume \eqref{TheCondition}. The the positive maximizer 
 $m^* = c^* >0$ of $\frg_\lambda^\beta $  is stable in the following sense: For every
 $c\in [0,\infty)$, 
\begin{equation}
 \label{gStability}
\frg_\lambda^\beta  (c^* ) - \frg_\lambda^\beta  (c)\, \geq\, \frac{\chi\beta (c^*)^3}{24} (c-c^* )^2 \, 
\equiv\, d(c - c^* ) .
\end{equation}
\end{lemma}

\section{Conclusion of the Proof of Theorem \ref{Main}: Claims \eqref{eq:stable} and \eqref{eq:decay}}
\label{sec:FullStability}
Let us return  first to the infinite dimensional variational problem \eqref{eq:dualvar}.
The FKG property of $\mu_\lambda^\beta$ (inherited in a limit of ferromagnetic 
Ising models) evidently implies that  $\Lambda_\lambda^\beta (h) \leq 
\Lambda_\lambda^\beta  (|h|)$.
Since by \eqref{Dorlas}, 
\[
 \Lambda_\lambda^\beta (|h|)\, -\, \frac1{2}\|h\|_\beta^2\, \leq\, 
\frac{1}{\beta}\int_0^\beta \frg_\lambda^\beta  (| h(t)|)\textd t, 
\]
we readily recover from \eqref{gStability} the following $L_2$-stability bound:
\begin{equation}
 \label{DualStability}
\begin{split}
\lbr \Lambda_\lambda^\beta (m^*\!\!\cdot\!\mathbf 1) - \frac{\beta (m^*)^2}2\rbr \, &- 
\, \lbr \Lambda_\lambda^\beta (|h|)\, -\, \frac1{2}\|h\|_{2;\beta}^2\rbr\, 
\\& \geq\, 
\frac1{\beta}\int_0^\beta d (|h (t)| - m^*)\textd t \, \equiv\, D (h) . 
\end{split}
\end{equation}
In order to transfer this bound to the original variational problem \eqref{eq:var} we 
shall again rely on duality considerations of Section~\ref{S:R-P}.  Recall that 
\begin{equation}
\Lambda_\lambda^\beta (m^*\!\!\cdot\! \mathbf 1 ) - \beta (m^* )^2/2 = \beta (m^* )^2/2 - 
I_{\lambda}^\beta (m^*\!\!\cdot \!\mathbf 1 ).
\end{equation}
Therefore, by \eqref{DualStability}
\[
 \beta (m^* )^2/2 - 
I_{\lambda}^\beta (m^*\!\!\cdot \!\mathbf 1 ) \, \geq\, 
\lbr (h ,m )_\beta - \| h\|^2_{
%2;
\beta}/2 \rbr  - I_{\lambda}^\beta  (m) + D(h) .
\]
Taking $h=m$  we arrive at
\begin{equation}
\label{DirStability}
\frG_{\lambda}^\beta ( m^*\!\!\cdot\!\mathbf 1 ) - \frG_{\lambda}^\beta (m)\, \geq\, D (m) .
\end{equation}
Obviously, there exists $c_1 =c_1 (\lambda ,\beta )$, such that 
\begin{equation}
 \label{Dbound}
 D (m)\, \geq\, c_1 \| m\pm m^*\cdot \mathbf 1\|_{
%2;
\beta}^2 .
\end{equation}
Thereby, the first term in \eqref{eq:stable} is recovered.  In order to recover the second
term in \eqref{eq:stable}, set $\bar m = m^*\!\!\cdot\! \mathbf 1$ and  just compute:
\begin{equation*}
 \begin{split}
  \frG_{\lambda}^\beta
 (\bar m) - \frG_{\lambda}^\beta
 (m)\, &=\, \lbr \frac12\|\bar m\|_{
%2;
\beta}^2 - \frac12\| m\|_{
%2;
\beta}^2\rbr
 - \lbr I^\beta_{\lambda} (\bar m ) - I^\beta_{\lambda} (m)\rbr\\
&=\, \frac12\|\bar m\|_{
%2;
\beta}^2 - \frac12\| m\|_{
%2;
\beta}^2 + (\bar m, m-\bar m )_\beta\\
&+ I^\beta_{\lambda } (m ) -  I^\beta_{\lambda} (\bar m )  - (\bar m, m-\bar m )_\beta\\
&=\, -\frac12\| m-\bar m\|_{
%2;
\beta}^2 \, +\, I^{\beta}_{\lambda ,m^*} (m) .
 \end{split}
\end{equation*}
Repeating the same computation with $-\bar m$, we infer:
\begin{equation*}
\begin{split} 
\frG_{\lambda}^\beta (\bar m) -
 \frG _{\lambda}^\beta(m)
\, &\geq \, - \frac12 \|m\pm \bar m\|_{
%2;
\beta}^2 + \min\lbr I^\beta_{\lambda ,m^*} (m) ,
I^\beta_{\lambda ,- m^*} (m) \rbr\, \\
&\geq\, - \frac12 \|m\pm \bar m\|_{
%2;
\beta}^2\, +\, 
\int_0^\beta U_{\eta(\lambda, m^*)} (m^\prime (t) )\textd t ,
\end{split}
\end{equation*}
as follows from \eqref{IStability}. At this point we recall 
\eqref{Dbound}, and the full statement of \eqref{eq:stable}  
becomes an easy exercise (with $U = U_\eta$).

Let now $m$ be such that $\| m\pm\bar m \|_{\rm sup}= 2\delta$. Let us consider only 
the case when $\max_t m(t) = m^* +2\delta$.  If $\min_t m(t) >\delta$, then
$\| m\pm\bar m\|_\beta^2\geq \delta^2\beta$. Otherwise, without loss of generality, 
we may assume that $m (0) = m^* +2\delta$ and let $r = \min\lbr t>0 : m(t) = m^* +\delta\rbr$.
Obviously, in such circumstances, 
\[
 \| m\pm\bar m\|_\beta^2 \geq r\delta^2\quad\text{whereas}\quad
\int_0^\beta U(m^\prime (t))\textd t\geq \int_0^rU(m^\prime (t))\textd t - 2\eta\beta .
\]
However, by convexity, 
\[
\int_0^rU(m^\prime (t))\textd t \geq r U (\delta /r )\sim \delta \log (1/r )  
\]
for $r\ll \delta$. Therefore, the $\int_0^\beta U(m^\prime (t))\textd t$-term in 
$\max\lbr \| m\pm\bar m\|_\beta^2  ,\int_0^\beta U(m^\prime (t))\textd t\rbr$ becomes 
dominant for $r\sim {\rm e}^{-1/\delta}$. Hence \eqref{supbound}.

Finally, consider the claim \eqref{eq:decay}.  
%Via the duality
As we have seen   in Section \ref{S:R-P}, $m^*!\!\cdot\!\mathbf 1$
% where $h^*$ 
solves the dual variational problem
\eqref{eq:dualvar}.  In particular,
\begin{equation}
m^*= \frac{1}{\beta} \mu^\beta_{\lambda, m^*}\left( \left( \sigma, \bf 1\right)_{\beta}\right)
\end{equation}
Via the Taylor expansion we may write
\begin{equation}
m^*= \frac{m^*}{\beta} \textrm{$\mathbb V$ar}^\beta_{\lambda, 0}\left( \left( \sigma, \bf 1\right)_{\beta}\right) - \frac{s_{4}\left(\lambda, \beta\right)(m^*)^3 }{6\beta} + o\left((m^*)^4\right)
\end{equation}
noting that odd derivatives of $\Lambda_{\lambda}^\beta $ vanish at zero by spin flip symmetry.
As we have already mentioned, 
 Lemma \ref{RCB} implies  that $s_4(\lambda, \beta) >0$.  Claim \eqref{eq:decay} follows easily.

\qed

%\begin{comment}
%\section{Quantum FK model on complete graphs}
%\label{FK}
%\end{comment}
\section{Appendix}

\begin{proof}[Proof of Lemma \ref{RCB}]
We begin by recalling the random current representation of the third
Ursell function for a discrete Ising model, see \cite{AF}.

Following the notation of \cite{AF}, consider a general finite graph
Ising model on a vertex set $V$ (also called sites), elements of which we denote with $i$, and with bonds $b$ in some fixed subset $\EE \subset V
\times V$. Let us denote the coupling constants by $J_b$ and local
fields by $h_i$.  It is useful to interpret  the local fields $h_i$ as the
coupling constants between the site $i$ and a `ghost' site
$\mathfrak g$.  We augment the graph $V, \EE$ by including $\frak g$ along with an edge set connecting $\frak g$ to each vertex in $V$.  Let us denote the augmented edge set by $\EE'$.

Let $\underline{n}= (n_b)_{b \in \EE'}$ denote a
sequence of integer valued `fluxes' attached to the bonds of the
underlying graph.  We say that that a site $i\in V$ is a \textit{source}
if $\sum_{i \in b} n_b$ is odd and denote the collection of sources
other than the ghost site by $\partial \underline{n}$.

We say that $x \leftrightarrow y$ if there exists a path of non-zero
fluxes connecting $x$ to $y$ using bonds in $\EE$ (i.e. not ghost
site bonds). Moreover, $x \nrightarrow h$ means that whenever
$x\leftrightarrow y$, $n_{y,\mathfrak g} = 0$. For any subset of
bonds $\mathcal B \subset \EE$, let
\begin{equation}
\bigl \langle \cdot \bigr \rangle_{\mathcal B}
\end{equation}
denote the Gibbs state for the Ising model with the constants
$\left(J_b\right)_{b\in \BB}$ set to $0$. Let
\begin{equation}
W(\underline{n})= \prod_{b} \frac{{J_b}^{n_b}}{n_b!}
\end{equation}
and
\begin{equation}
Z= \sum_{\partial \underline{n}= \varnothing} W(\underline{n}).
\end{equation}

Suppose $r, s, t \in V$ are fixed sites in the Ising model under
consideration.  To avoid excessive provisos hereafter we shall assume that all three points are distinct.  For a collection of fluxes $\underline{n}$, let us
say that $b \in C_{\underline n}(r)$ if $i \leftrightarrow r$ for
some $i \in b$.  Then we have (the reader should consult \cite{AF} for a derivation)
\begin{multline}\label{3RCR}
U_{\lambda, c}^\beta(r,s,t) = \\
\left\{\sum_{\partial \underline{n}_1= \{r, s\}, \:
\partial \underline{n}_2 = \varnothing}
\frac{W(\underline{n}_1)}{Z} \frac{W(\underline{n}_2)}{Z} \mathbf
1_{\underline{n}_1 + \underline{n}_2 : r \nrightarrow h} \times
\left[\bigl \langle \sigma_t \bigr \rangle_{C^c_{\underline{n}_1 +
\underline{n}_2}(r)} - \bigl \langle \sigma_t \bigr \rangle\right]
\right\} + \left\{s \Leftrightarrow t\right\}
\end{multline}
where $\left\{s \Leftrightarrow t\right\}$ represents the first term
with the roles of $s$ and $t$ interchanged.  Here $\underline{n}_1$
and $\underline{n}_2$ are independent copies of fluxes.

Clearly, the weights $W(\underline{n})$ are proportional to the
probability that a family of independent Poisson processes indexed
by (generalized) bonds take a collection of values determined by
$\underline n$. We need to differentiate between the processes
associated to the bonds of the graph and the bonds with the ghost
site $\mathfrak g$. Specifically, let $\{\mathcal N_b, \mathcal
M_i\}_{b \in \EE, i \in V}$ denote a collection of independent
Poisson processes with respective parameters \newline $\{J_b,
h_i\}_{b \in \EE, i \in V}$ and let $\textd \mathbb P$ denote the joint
probability measure associated to these processes.

It is well known that if $h \geq 0$, spin correlations are
increasing with respect to coupling strengths, so each summand on
the righthand side of \eqref{3RCR} must be non-positive. Therefore, neglecting
summands we obtain
\begin{multline}\label{rcbound1}
U_{\lambda, c}^\beta(r,s,t) \leq \\
\sum_{\partial \underline{n}_1\{r, s\}, \:
\partial \underline{n}_2 = \varnothing }
\frac{W(\underline{n}_1)}{Z} \frac{W(\underline{n}_2)}{Z} \mathbf
1_{\underline{n}_1 + \underline{n}_2 : r \nrightarrow h}
1_{\underline{n}_1 + \underline{n}_2 : r \leftrightarrow t} \times
\left[\bigl \langle \sigma_t \bigr \rangle_{C^c_{\underline{n}_1 +
\underline{n}_2}(r)} - \bigl \langle \sigma_t \bigr \rangle\right] +
\left\{s \Leftrightarrow t\right\}
\end{multline}
where we drop analogous terms in the expression $\left\{s
\Leftrightarrow t\right\}$.

Let us concentrate on the first term on the righthand side of the
inequality as the second may be treated by symmetric considerations.  We note two things. First, on the set $\bigl
\{\underline{n}_1 + \underline{n}_2 : r \leftrightarrow t \bigr\}
\cap \bigl \{\underline{n}_1 + \underline{n}_2 : r \nrightarrow h
\bigr \}$ we have
\begin{equation}
\bigl \langle \sigma_t \bigr \rangle_{C^c_{\underline{n}_1 +
\underline{n}_2}(r)} = 0
\end{equation}
by spin flip symmetry. Second, we have
\begin{multline}\label{requirement}
\bigl\{\partial \underline{n}_1= \{r, s\}, \partial \underline{n}_2
= \varnothing \bigr\}\cap \bigl\{\underline{n}_1 + \underline{n}_2 :
r
\nrightarrow h\bigr \} = \\
\bigl\{\underline{n}_1 : r \leftrightarrow s \bigr\} \cap
\bigl\{\partial \underline{n}_1= \{r, s\},
\partial \underline{n}_2 = \varnothing \bigr\} \cap\bigl\{\underline{n}_1 + \underline{n}_2 : r \nrightarrow h\bigr\}.
\end{multline}

In terms of the Poisson processes, let $E_{r,s,t}$ denote the event
determined by the requirements of the first term on the righthand
side of \eqref{rcbound1}.  After the reduction made above,
\eqref{rcbound1} may be expressed as
\begin{equation}
\label{EQ:UB}
U_{\lambda, c}^\beta(r,s,t) \leq -\langle\sigma_t \rangle  \frac{\mathbb P\otimes \mathbb P (E_{r,s,t})}
{\left(\mathbb P (\underline{n}_1 = \varnothing)\right)^2} + \left\{s
\Leftrightarrow t\right\}
\end{equation}

Let us now specialize to Ising models on the circle with $2^N$
vertices.  We take coupling and field strengths given by
\begin{align*}
\texte^{-2J_b} & = \frac{\lambda \beta}{2^N}\\
h_i &= \frac{h\beta}{2^N}.
\end{align*}
For the readers convenience we note that
\begin{align*}
\mathbb P( \mathcal N_b = \text{odd}) & \propto \tanh(J_b) =
\frac{1- \frac{\lambda \beta}{2^N}}{1+\frac{\lambda \beta}{2^N}}\\
\mathbb P( \mathcal N_b = \text{even}) & \propto 1 \\
\mathbb P( \mathcal N_b = 0) \propto  & \cosh(J_b)^{-1} =
\sqrt{\frac{\lambda \beta}{2^N}} + O(2^{-N})
\end{align*}
and
\begin{align*}
\mathbb P( \mathcal M_i = \text{odd}) & = \frac{h \beta}{2^N}e^{-\frac{h \beta}{2^N}} + O(2^{-3N})\\
\mathbb P( \mathcal M_i = \text{even}) & =  e^{-\frac{h \beta}{2^N}} + O(2^{-2N}) \\
\mathbb P( \mathcal M_i = 0)  & = e^{-\frac{h \beta}{2^N}}
\end{align*}
where the implicit constant of proportionality is
\begin{equation}
\label{eq:const-prop}
(1+ o(1))e^{-\ffrac{\lambda \beta}{2^N}}.
\end{equation}

We compute a lower bound for the probability determined by the
numerator in \eqref{EQ:UB}. 
Without loss of generality we may assume the
orientation
\begin{equation}
0<r <s < t < \beta.
\end{equation}
Continuity arguments imply that it is sufficient for us prove the
bound \eqref{thebound} assuming each point is of the form $\ffrac{ \beta j}{ 2^{k}}$ for some $j, k\in \mathbb N$ fixed.

Let $c, d \in [0, \beta]$ be fixed.
Let $\lfloor c, d \rfloor$ denote the
integers $i$ in $[0,2^N]$ such that $\frac{c 2^N}{\beta} \leq i
\leq \frac{d 2^N}{\beta}$. If $c>d$ this is understood to be with
respect to periodic boundary conditions.  For any bond $b$ let $b\in
\lfloor c, d \rfloor$ denote the case that both vertices of $b$ are
in $\lfloor c, d \rfloor$. Finally, let
\begin{equation}
M_{\lfloor c, d \rfloor}= \{\mathcal M^{(1)}_i + \mathcal M^{(2)}_i
> 0 \text{ for some $i \in \lfloor c, d \rfloor$}\}.
\end{equation}

Based on topological considerations there are three cases: $E_{r,s,t}\cap M_{\lfloor r, s\rfloor},E_{r,s,t}\cap M_{\lfloor s, t\rfloor}$ and
$E_{r,s,t}\cap M_{\lfloor t, r\rfloor}$. The requirements of $E_{r,s,t}$
imply these three possibilities lead to disjoint subevents (but not a partition) of
$E_{r,s,t}$. 

The cases carry a minimal but tedious amount of
computation, so we shall present only one of them.

Let us consider $E_{r,s,t} \cap M_{\lfloor t, r \rfloor}$.  For any $x <r<
t< y$ in the dyadic lattice $\{\frac {\beta l}{2^N}\}$ , let
\begin{multline}
A(x,y)= \{ \mathcal N^{(1)}_b + \mathcal N^{(2)}_b > 0 \text{ for
all $b \in \lfloor x , r \rfloor \cup \lfloor t, y \rfloor$}\}\cap \\
\{\mathcal N^{(1)}_b + \mathcal N^{(2)}_b = 0, \:\text{ if } b = \{\ffrac{(x-1)
2^N}{\beta},  \ffrac{x 2^N}{\beta}\} \text{ or } \{\ffrac{y
2^N}{\beta}, \ffrac{(y+1)2^N}{\beta}\}\}.
\end{multline}

The parity of all b-bonds in $\lfloor x, y \rfloor$ is completely
determined on the event $E_{r,s,t} \cap M_{\lfloor t, r \rfloor}\cap
A(x,y).$  Since the fluxes associated to $\frak g$-bonds must be zero, independence of the various Poisson processes allows us to
compute:
\begin{equation}
\mathbb P\left(E_{r,s,t} \cap M_{\lfloor r, t \rfloor} \cap A(x,y) \right)
\propto  e^{-2\lambda |s-r| - 2 h |y-x|}
\left(\frac{\lambda \beta}{2^N} \right)^2 + o(2^{-2N}).
\end{equation}
A worst case lower bound for the exponential term is $e^{-(2\lambda+
2 h) \beta}$, so summing over ordered pairs $\{x, y\}$ and taking into account \eqref{eq:const-prop},
\begin{equation}
\mathbb P\left(E_{r,s,t} \cap M_{\lfloor r, t \rfloor}\right) \geq
e^{-(4\lambda + 2h)\beta} \frac{1}{2}[\lambda
\left(\beta+r-t\right)]^2 + o(1).
\end{equation}
Analogous estimates in the other cases lead to
\begin{equation}
\mathbb P\left(E_{r,s,t} \right) \geq \frac{\lambda^2}{2} e^{-(4
\lambda + 2 h) \beta} \left[ (s-r)^2 + (t-s)^2 + (\beta+r-t)^2\right] +
o(1).
\end{equation}
As $\mathbb P (\partial \underline{n} = \varnothing) \leq 1$, the
bound \eqref{thebound} follows.

\begin{comment}
\begin{multline}
\mathbb P \left(\mathcal N^{(1)}_b = \text{ odd for all $b \in
\lfloor r,s \rfloor$}, \: \mathcal N^{(1)}_b + \mathcal N^{(2)}_b >
0 \text{ for all $b \in \lfloor s,t \rfloor$}\right) \mathbb P
\left(\mathcal M^{(1)}_i + \mathcal M^{(2)}_i = 0 \text{ for all $i
\in \lfloor r,t \rfloor$}, \mathcal M^{(1)}_i + \mathcal M^{(2)}_i >
0 \text{ for some $i \in
\end{multline}
\end{comment}

\end{proof}


\begin{thebibliography}{aaa}

\bibitem{Aizenman}
M.~Aizenman, \textit{Geometric analysis of $\varphi^4$ fields and
Ising models. I, II.},  Commun. Math. Phys.  \textbf{86} (1982)
1--48.

\bibitem{ABF}
M.~Aizenman, D.J.~Barsky and R.~Fern\'andez, \textit{The phase
transition in a general class of Ising-type models is sharp},
J.~Statist.~Phys. \textbf{47} (1987) 343--374.


\bibitem{ACCN}
M.~Aizenman, J.T.~Chayes, L.~Chayes and C.M.~Newman,
\textit{Discontinuity of the magnetization in one-dimensional
\hbox{$1/\vert x-y\vert^2$} Ising and Potts models}, J.~Statist.
Phys. \textbf{50} (1988), no. 1-2, 1--40.

\bibitem{AF}
M.~Aizenman and R.~Fern\'andez, \textit{On the critical behavior of
the magnetization in high-dimensional Ising models}  J.~Statist.
Phys.  \textbf{44} (1986) 393--454.

\bibitem{AKN}
M. Aizenman, A. Klein, C. Newman,  \textit{Percolation methods for
disordered quantum Ising models}, in Phase Transitions: Mathematics,
Physics, Biology,..., R. Kotecky,ed., pp. 1-26, World Scientific,
Singapore 1993.

\bibitem{AN}
M.~Aizenman, B.~Nachtergaele, \textit{Geometric aspects of quantum
spin states}, Commun. Math. Phys.   164 , 17--63 (1994)

\bibitem{Baldi} P.~Baldi, \textit{
Large deviations and stochastic homogenization}, 
Ann. Mat. Pura Appl. (4) 151, 161--177 (1988).

%\bibitem{Bill}
%P.~Billingsley, \textit{Probability and Measure} New York, NY: John
%Wiley \& Sons (1995)

%\bibitem{B-K}
%M. Biskup and R. Koteck\'{y}, \textit{Forbidden gap argument for
%phase transitions proved by means of chessboard estimates}, Commun.
%Math. Phys.  264  (2006), no. 3, 631-656.

\bibitem{Cam-Kl-Per}
M. Campanino, A.~Klein, J.F. Perez \textit{Localization in the
ground state of the Ising model with a random transverse
field}, Communications in Mathematical Physics 135, 499-515 (1991).

%\bibitem{C-H}
%Carmona, P.; Hu, Yueyun \textit{Universality in Sherrington
%Kirkpatrick's spin glass model}, Annales de l'Institut Henri
%Poincar\'{e} (B) Probability and Statistics, \textbf{42}, (2006),
%no. 2, 215--222,  Elsevier

%\bibitem{CFMP}
%M.~Cassandro, P.A.~Ferrari, I.~Merola and E. Presutti,
%\textit{Geometry of contours and Peierls estimates in $d=1$ Ising
%models}, math-ph/0211062.

%\bibitem{Cassandro-Presutti}
%M.~Cassandro and E.~Presutti, \textit{Phase transitions in Ising
%systems with long but finite range interactions}, Markov Process.
%Related Fields \textbf{2} (1996) 241--262.

%\bibitem{CKS1}
%L. Chayes, R. Koteck\'y and S. B. Shlosman, \textit{Aggregation and
%intermediate phases in dilute spin systems}, Communications in
%Mathematical Physics \textbf{171} (1995), no. 1, 203-232.

%\bibitem{CKS2}
%L. Chayes, R. Koteck\'{y}, S. Shlosman,  \textit{Staggered Phases in
%Diluted Systems with Continuous Spins}, Commun. Math. Phys.
%\textbf{189}, 631-640 (1997)

%\bibitem{Dembo-Zeit}
%A. Dembo and O. Zeitouni \textit{Large deviations techniques and applications} Springer, New York, (1998).

%\bibitem{Dobrushin-a}
%R.~Doburshin, \textit{Gibbsian random fields for lattice systems
%with pair interaction.} Funct. Anal. and Appl., 2 (1968), N4, 31-43.

%\bibitem{Dobrushin-b}
%R.~Dobrushin, \textit{The problem of uniqueness of a Gibbs random
%fields and the problem of phase transition.} Func. anal. and appl.,
%2, p. 302-312, 1968.

\bibitem{Dorlas}
T.C.~Dorlas, \textit{Probabilistic derivation of a noncommutative
version of Varadhan's theorem}, unpublished, June 2002,
\verb"http://www.stp.dias.ie/~dorlas/tony_index2.html"

%\bibitem{Ellis-at-al}
%M.~Costeniuc, R.S.~Ellis and H.~Touchette, \textit{Complete analysis
%of phase transitions and ensemble equivalence for the
%Curie-Weiss-Potts model}, preprint 2004
%(http://xxx.lanl.gov/abs/cond-mat/0410744).

%\bibitem{Cug-Grem}
%L.~Cugliandolo, D.~Grempel, C.~da Silva Santos,
%\textit{Imaginary-time replica formalism study of a quantum
%spherical p-spin-glass model} Phys. Rev. B \textbf{64}, 014403
%(2001).

%\bibitem{DS1}
%R. L. Dobrushin and S. B. Shlosman, \textit{Absence of Breakdown of
%Continuous Symmetry in Two-dimensional Models of Statistical
%Physics}, Commun. Math. Phys. \textbf{42}, 31, (1975)

%\bibitem{DS}
%R.L.Dobrushin and S.B. Shlosman, \textit{Phases corresponding to the
%local minima of the energy}, Sel. Math. Sov. 1, 31-338, 1981.

%\bibitem{Dor}
%T.C. Dorlas, \textit{Probabilistic derivation of a noncommutative
%version of Varadhan's theorem}, June 2002 preprint

%\bibitem{Dyson1}
%F.J.~Dyson, \textit{Existence of a phase-transition in a
%one-dimensional Ising ferromagnet}, Commun. Math. Phys. \textbf{12}
%(1969), no. 2, 91--107.

%\bibitem{Dyson2}
%F.J.~Dyson, \textit{An Ising ferromagnet with discontinuous
%long-range order}, Commun. Math. Phys. \textbf{21} (1971), 269--283.

%\bibitem{DLS}
%F.J.~Dyson, E.~Lieb and B.~Simon, \textit{Phase Transitions in
%Quantum Spin Systems with Isotropic and Non-Isotropic Interactions},
%Journal of Statistical Physics, \textbf{18}, 335--383 (1978).

%\bibitem{Ellis}
%R. S. Ellis \textit{Entropy, Large Deviations, and Statistical
%Mechanics}, Grundlehren der mathematischen Wissenschaften 271,
%Springer-Verlag, New York Berlin Heidelberg Tokyo, 1985

%\bibitem{vESh}
%A. C. D. van Enter and S. B. Shlosman, \textit{Provable first-order
%transitions for nonlinear vector and gauge models with continuous
%symmetries}, Communications in Mathematical Physics 255 (2005), no.
%1, 21-32.

\bibitem{FSW}
M. Fannes, H. Spohn, and A. Verbeure, \textit{Equilibrium states for mean field 
models}, Journal of Mathematical Physics 21, 355�360 (1980)

\bibitem{FFS}
R.~Fern\'andez, J.~Fr\"ohlich and A.D~Sokal, \textit{Random walks,
critical phenomena, and triviality in quantum field theory}, Texts
and Monographs in Physics, Springer-Verlag, Berlin, 1992.

%\bibitem{Froh}
%J.~Fr{\"o}hlich \textit{On the triviality of $d > 4$ theories and
%the approach to the critical point in d(-) > 4 dimensions}, Nuclear
%Physics B, Volume 200, Issue 2, 1 February 1982, Pages 281-296
%1981. Available online 21 October 2002


\bibitem{FILS1}
J.~Fr{\"o}hlich, R.~Israel, E.H.~Lieb and B.~Simon, \textit{Phase
transitions and reflection positivity. I.~General theory and
long-range lattice models}, Commun. Math. Phys. \textbf{62} (1978),
no. 1, 1--34.

\bibitem{FILS2}
J.~Fr{\"o}hlich, R.~Israel, E.H.~Lieb and B.~Simon, \textit{Phase
transitions and reflection positivity. II.~Lattice systems with
short-range and Coulomb interations}, J.~Statist. Phys. \textbf{22}
(1980), no. 3, 297--347.

%\bibitem{FL}
%J.~Fr\"{o}hlich and E.H.~Lieb, \textit{Ph036530285  דינא דאהרase transitions in
%anisotropic lattice spin systems}, Commun. Math. Phys.,\textbf{60},
%3, (1978), 233--267

%\bibitem{FSS}
%J.~Fr{\"o}hlich, B.~Simon and T.~Spencer, \textit{Infrared bounds,
%phase transitions and continuous symmetry breaking}, Commun. Math.
%Phys. \textbf{50} (1976)   79--95.

%\bibitem{Georgii}
%H.-O.~Georgii, \textit{Gibbs Measures and Phase Transitions},
%de~Gruyter Studies in Mathematics, vol.~9, Walter de Gruyter~\&~Co.,
%Berlin, 1988.

\bibitem{G-P-S}
A.~Georges, O.~Parcollet, S.~Sachdev, \textit{Mean field theory of a
quantum Heisenberg spin glass} Physical Review Letters \textbf{85},
840 (2000)

\bibitem{Gin}
J.~Ginibre \textit{Existence of phase transitions for quantum
lattice system} Commun. Math. Phys. \textbf{14} (1969), 205.

\bibitem{Grim}
G.~Grimmett, \textit{Space-time percolation} Preprint,
arXiv:0705.0506v1 [math.PR]



%\bibitem{G-J}
%J. Glimm and A. Jaffe, \textit{Quantum Phy036530285  דינא דאהרsics}, Springer-Verlag,
%New York, 415 pages (1981).

%\bibitem{G-J-S}
%J. Glimm, A. Jaffe, T. Spencer, \textit{Phase Transitions for
%$\phi^4_2$ Quantum Fields}, Commun. Math. Phys., 45 (1975), 203�216.

%\bibitem{Grif}
%R. Griffiths, Phys. Rev. 152:240-246 (1966).

\bibitem{Guerra}
F.~Guerra, \textit{Broken replica symmetry bounds in the mean field
spin glass model}, Commun. Math. Phys. \textbf{233} (2003), no. 1,
1--12.

%\bibitem{Guerra-Ton-1}
%F.~Guerra, F.L.~Toninelli, \textit{The thermodynamic limit in mean
%field spin glass models}, Commun. Math. Phys. \textbf{230} (2002) 1,
%71-79

%\bibitem{Guerra-Ton-2}
%F.~Guerra,  F.L.~Toninelli, \textit{The infinite volume limit in
%generalized mean field disordered models} Markov Process. Related
%Fields  \textbf{9} (2003), no. 2, 195--207

%\bibitem{Haldane}
%F. D. M. Haldane, \textit{Continuum dynamics of the 1-D Heisenberg
%antiferromagnet: Identification with the $O(3)$ nonlinear sigma
%model} Phys. Lett. 93A , 464--468 (1983).

%\bibitem{Hara-Remco-Slade2}
%T. Hara, R. van der Hofstad, and G. Slade. \textit{Critical
%two-point functions and the lace expansion for spread-out
%high-dimensional percolation and related models}, Ann. Probab.
%\textbf{31} (2003), no.~1, 349--408.

%\bibitem{Hara-Slade1} T.~Hara and G.~Slade,
%\textit{Self-avoiding walk in five or more dimensions. I.~The
%critical behaviour}, Commun. Math. Phys.~\textbf{147} (1992)
%101--136.

%\bibitem{HL}
%O.~Heilmann and E.~Lieb, \textit{Lattice models for liquid
%crystals}, Journal of Statistical Physics, \textbf{20}, 6, 679--693,
%(1979)

\bibitem{ILN} D.~Ioffe, \textit{Stochastic geometry of classical 
and quantum Ising models}, preprint (2007).

\bibitem{Ioffe}
D.~Ioffe and   A.~Levit, \textit{  Long range order and giant components
of quantum random graphs},  Mark. Proc. Rel. Fields \textbf{13} (2007), no. 3,
469--492.

%\bibitem{Remco-Gordon} R.~van~der~Hofstad and G.~Slade,
%\textit{A generalised inductive approach to the lace expansion},
%Probab. Theory Rel. Fields~\textbf{122} (2002) 389--430.

%\bibitem{Remco-Gord2}
%R.~van der Hofstad and G.~Slade, \textit{Convergence  of critical
%oriented percolation to super-Brownian motion above $4+1$
%dimensions}, Ann. Inst. H. Poincar� Probab. Statist. \textbf{39}
%(2003), no. 3, 413--485.

%\bibitem{Remco-Frank-Gord}
%R.~van~der~Hofstad, F.~den Hollander, G.~Slade, \textit{Construction
%of the incipient infinite cluster for spread-out oriented
%percolation above $4+1$ dimensions}, Commun. Math. Phys.
%\textbf{231}  (2002), no. 3, 435--461.

%\bibitem{Kenn}
%T. Kennedy, \textit{Long range order in the anisotropic quantum
%ferromagnetic Heisenberg model},  Commun. Math. Phys. \textbf{100},
%447 (1985)

%\bibitem{Ken-Tas}
%Kennedy, T. and Tasaki, H. \textit{Hidden symmetry breaking and the
%Haldane phase in s = 1 quantum spin chains}, Comm. Math. Phys. 147
%(1992), 431�484.

%\bibitem{Kerimov93}
%A.~Kerimov, \textit{Absence of phase transitions in one-dimensional
%antiferromagnetic models with long-range interactions},  J.~Statist.
%Phys. \textbf{72} (1993), no. 3-4, 571--620.

%\bibitem{Kerimov96}
%A.~Kerimov, \textit{Phase transition in one-dimensional model with
%unique ground state}, Physica~A \textbf{225} (1996), no. 2,
%271--276.

\bibitem{MPV}, preprint.

M.~M\'ezard,G.~Parisi, and M.A.~Virasoro, \textit{Spin glass theory
and beyond}, World Scientific Lecture Notes in Physics, vol.~9,
World Scientific Publishing Co., Inc., Teaneck, NJ, 1987.

\bibitem{N-1}
B. Nachtergaele, \textit{Quasi-state decompositions for quantum spin
systems in  Probability Theory and Mathematical Statistics}
(Proceedings of the 6th Vilnius Conference) , B. Grigelionis et al.
(Eds), VSP/TEV, Utrecht-Tokyo-Vilnius, 1994, pp 565-590

\bibitem{N-2}
B.~Nachtergaele, \textit{A stochastic geometric approach to quantum
spin systems}.  Probability and phase transition (Cambridge, 1993),
237--246, NATO Adv. Sci. Inst. Ser. C Math. Phys. Sci., 420, Kluwer
Acad. Publ., Dordrecht, 1994.

\bibitem{Parisi}
G.~Parisi, \textit{Field theory, disorder and simulations}, World
Scientific Lecture Notes in Physics, \textbf{49}, World Scientific
Publishing Co., Inc., River Edge, NJ, 1992.

\bibitem{S-K}
D.Sherrington and S.Kirkpatrick, \textit{Solvable model of a
spin-glass}, Phys. Rev. Lett., 35, 1792-1796 (1975)

\bibitem{Talagrand-book}
M.~Talagrand, \textit{Spin Glasses: A Challenge for Mathematicians.
Cavity and Mean Field Models}, A Series of Modern Surveys in
Mathematics, vol~46. Springer-Verlag, Berlin, 2003.

\bibitem{Talagrand-paper}
M.~Talagrand, \textit{The Parisi formula}, Ann. of Math. (2)
\textbf{163} (2006), no. 1, 221--263.

\bibitem{Toland} J.F.~Toland, \textit{A duality pronciple for 
non-convex optimization in the calculus of variations}, F.M.R.I. (University of Essex), Arch. Rational Mech. Analysis, 1979. 
\end{thebibliography}
\end{document}